\documentclass[a4paper, 11pt]{amsart}
\usepackage{a4wide}
\usepackage{hyperref}
\usepackage{graphicx}
\usepackage{amsmath,amssymb,amsfonts}
\usepackage{amsthm}
\usepackage{mathrsfs}
\usepackage{xcolor}
\usepackage{booktabs}
\usepackage{algorithm}
\usepackage{algpseudocode}
\algrenewcommand\algorithmicrequire{\textbf{Input:}}
\algrenewcommand\algorithmicensure{\textbf{Output:}}

\newtheorem{prop}{Proposition}[section]
\newtheorem{corollary}[prop]{Corollary}

\newtheorem{lemma}[prop]{Lemma}
\newtheorem{thm}[prop]{Theorem}

\theoremstyle{definition}
\newtheorem{rem}[prop]{Remark}
\newtheorem{ques}[prop]{Question}

\newcommand{\Aut} {\mathrm{Aut}}
\newcommand{\Isom} {\mathrm{Isom}}
\newcommand{\End} {\mathrm{End}}

\newcommand{\ord}{\mathrm{ord}}
\newcommand{\nrd}{\mathrm{nrd}}
\newcommand{\eps}{\varepsilon}
\def\BF{\mathbb{F}}
\def\BZ{\mathbb{Z}}

\def\CO{\mathcal{O}}
\def\Fp{\mathfrak{p}}

\title{Divisibility of the coefficients of Drinfeld modular polynomials}

\author{Florian Breuer \and Fabien Pazuki}

\address{School of Computer and Information Sciences, The University of Newcastle,
University Drive, Callaghan, NSW 2308, Australia.}
 \email{Florian.Breuer@newcastle.edu.au}

\address{Department of Mathematical Sciences, University of Copenhagen,
 Universitetsparken 5, 
2100 Copenhagen \O, Denmark, and Universit\'e de Bordeaux, 33405 Talence, France.}
 \email{fpazuki@math.ku.dk}

\begin{document}

\subjclass[2020]{11G09; 11F52}
\keywords{Drinfeld modules, Drinfeld modular polynomials, isogenies, complex multiplication, singular moduli}

\begin{abstract}
    Let $A = \BF_q[T]$ and let $N \in A$ be monic of positive degree. Let $\Phi_N(X,Y)\in A[X,Y]$ be the Drinfeld modular polynomial which vanishes precisely at pairs of $j$-invariants of rank~2 Drinfeld $A$-modules linked by a cyclic isogeny of degree $N$. We study the divisibility of the coefficients of $\Phi_N(X+J_1, Y+J_2)$ by primes $P \in A$ for certain $J_1$, $J_2$, in particular $J_1=J_2=0$. These coefficients turn out to be highly divisible by primes $P$ of small degree at which the Drinfeld modules associated to $J_1, J_2$ have supersingular reduction. We also present a new algorithm which computes $\Phi_N$ exactly, by deforming the Drinfeld module with $j$-invariant $0$ at a single place of $A$, the same deformations used for our proofs. Finally we use this algorithm to test our results numerically.

\end{abstract}

\maketitle


\section{Introduction} 
Let $q$ be a power of a prime $p$, let $A = \BF_q[T]$ and $k = \BF_q(T)$. Let $k_\infty = \BF_q((1/T))$ denote the completion of $k$ at the infinite place and $\mathbb{C}_\infty$ the completion of an algebraic closure of $k_\infty$.

A rank 2 Drinfeld $A$-module $\varphi$ over an $A$-field $L$ (i.e.\ a field $L$ equipped with an $A$-algebra structure $\gamma : A \to L$) is determined by
\[
\varphi_T = \gamma(T) + g\tau + \Delta\tau^2, \quad \mathrm{where}\quad g, \Delta \in L, \; \Delta \neq 0,
\]
and where $\tau$ denotes the $q$-power Frobenius $x \mapsto x^q$. Its $j$-invariant is defined by $$j(\varphi) = \frac{g^{q+1}}{\Delta},$$ and it determines $\varphi$ up to isomorphism over an algebraic closure of $L$. For standard references on Drinfeld modules, we refer to \cite{Goss1996, Papikian2023}.

Let $N\in A$ be monic of positive degree. An isogeny $u : \varphi \to \varphi'$ between two Drinfeld modules of rank 2 is called \emph{cyclic of degree $N$} if its kernel is isomorphic to $A/NA$ as an $A$-module. There exists a \emph{Drinfeld modular polynomial}
\[
\Phi_N(X,Y) \in A[X,Y],
\]
symmetric and monic. Let us denote $|N| := q^{\deg N} = \#(A/NA)$. Then $\Phi_N$ has degree in $X$ (and hence also in $Y$) equal to
\[
\psi(N) = |N|\prod_{P \mid N}\big(1 + |P|^{-1}\big),
\]
and vanishes precisely at the pairs $(j(\varphi), j(\varphi'))$ of $j$-invariants of rank 2 Drinfeld modules linked by a cyclic $N$-isogeny; here $P$ runs over the monic prime divisors of $N$. The polynomial $\Phi_N$ was constructed by Bae \cite{Bae1992} (see also \cite{BR2009, BR2016} for arbitrary rank), and the case $N = T$ was studied in detail by Schweizer \cite{Schweizer1995} and, more recently, in \cite{BRaz2024}. As in the classical case, the coefficients of $\Phi_N$ grow rapidly: the height $h(\Phi_N)$, defined as the maximal $T$-degree of these coefficients, satisfies the following estimate:

    \begin{equation}\label{explicitbound}
    h(\Phi_N) = \frac{q^2-1}{2}\psi(N)\left(\deg N - 2\lambda_N + b_q(N)\right),
    \end{equation}
    where
    \[
    \lambda_N = \sum_{P^n\|N}\frac{|P|^n-1}{|P|^{n-1}(|P|^2-1)}\deg P
    \qquad\text{and}\qquad
    0\leq b_q(N)\leq B_q := 4 + \frac{5q^3-q^2-q+1}{q(q+1)(q-1)^2},
    \]
    see \cite[Thm.~1.1]{BPR2024}. The asymptotic form of \eqref{explicitbound}, with an inexplicit bounded term, is due to Hsia \cite{Hsia1998}; see also \cite{CGS2020, BRaz2024}, the latter computing the height exactly at level $N = T$.

\begin{ques}
     This height is huge, so at least one coefficient must grow fast with $\deg N$. Can one be more precise? For instance, is there a big prime factor explaining this size? Or does it come from a large number of small factors? 
\end{ques}

 This paper is concerned with the \emph{divisibility} of the coefficients of $\Phi_N$ by small primes of $A$. It is the Drinfeld module analogue of \cite{Breuer2025}, which in turn was motivated by conjectures and results of Wang \cite{Wang2026} on the classical modular polynomials: indeed in the elliptic curve case, modular polynomials are often shown with factorised coefficients, and they tend to have a lot of small factors. What is the situation in the case of Drinfeld modules? Our first main result is the following; we write $v_P$ for the normalized $P$-adic valuation on $k$.

\begin{thm}\label{thm:main}
    Let $N \in A$ be monic of positive degree and write
    \[
    \Phi_N(X,Y) = \sum_{0\leq m,n\leq \psi(N)} a_{m,n}X^mY^n, \qquad a_{m,n}\in A.
    \]
    Let $P\in A$ be a monic prime with $P \nmid N$. Then for all $m + n < \psi(N)$ the following hold.
    \begin{enumerate}
        \item If $\deg P = 1$, then
        \[
        v_P(a_{m,n}) \;\geq\; (q+1)\big(\psi(N) - m - n\big).
        \]
        \item If $\deg P$ is odd, then
        \[
        v_P(a_{m,n}) \;\geq\; (q+1)\big(C_0(N,P) - m - n\big),
        \; \mathrm{where}\; C_0(N,P) := \ord_X\big(\Phi_N(X,0)\bmod P\big).
        \]
    \end{enumerate}
\end{thm}

Part (1) is the special case of part (2) in which $C_0(N,P) = \psi(N)$; this holds for $\deg P = 1$ because $j=0$ is then the unique supersingular $j$-invariant in characteristic $P$ (Proposition \ref{prop:fullpsi}). For $\deg P$ even, the reduction of $j = 0$ is ordinary and no such divisibility occurs (Proposition \ref{ordinary}); for $\deg P$ odd, $C_0(N,P)$ is computed by the arithmetic of a quaternion algebra over $k$ (Proposition \ref{supersingular}(1)). In this case Theorem \ref{thm:main}(2) gives non-trivial divisibility of non-zero coefficients of $\Phi_N(X,Y)$ only when $C_0(N,P) > C_0(N,0)$ (Remark \ref{rem:vanishing}), and for odd $q$ this forces $\deg P\leq\deg N$ (Proposition \ref{supersingular}(2)); such excess does occur (Section \ref{sec:higherdeg}).

Theorem \ref{thm:main} is the analogue of \cite[Theorem~1.1]{Breuer2025}, where it is shown that the coefficients of the classical modular polynomial $\Phi_N^{\mathrm{class}}(X,Y)$ satisfy $v_p(a_{m,n}) \geq n_p(C_0(N,p
)-m-n)$ for $p\geq2$ prime, $p\nmid N$, and $n_p$ is an explicit integer. It is instructive to compare the two situations. In the elliptic case, the constants $n_p$ arise from the automorphisms of the elliptic curve $E$ with $j$-invariant $0$: one has $n_p = \tfrac{1}{2}\#\Aut(E_{\bar\BF_p})$ for $p \geq 5$, while for $p = 2$ or $p=3$ wild ramification phenomena increase $n_p$ further and introduce a dependence on $N$ modulo small integers. In the Drinfeld module setting, the automorphism group functor is always \emph{tame}: the rank 2 Drinfeld module $\phi^0$ with $j(\phi^0) = 0$ has automorphism group $\BF_{q^2}^*$, of order prime to the characteristic, and
\[
q + 1 = \frac{\#\BF_{q^2}^*}{\#\BF_q^*} = \frac{\#\Aut(\phi^0)}{\#\Aut(\phi^{\mathrm{gen}})}
\]
plays the role of $n_p$ from the classical case. Consequently the constant $q+1$ in Theorem \ref{thm:main} is uniform: it depends neither on $P$ nor on $N$. Computations suggest that it is also optimal for every $N$; see for instance Table \ref{tab:PhiT2T1} in Section \ref{sec:computations} to see an example in fully factored form.

More generally, we study the coefficients of $\Phi_N(X + J_1, Y + J_2)$ for $j$-invariants $J_1, J_2$ of Drinfeld modules with good reduction over local fields (Theorem \ref{thm:bigprimes}). The CM field of $j = 0$ is the constant field extension $\BF_{q^2}(T)$ of $k=\mathbb{F}_q(T)$, which is everywhere unramified, and this is the source of the uniformity in Theorem \ref{thm:main}. By contrast, at a degree one prime $P$ where a $j$-invariant $J$ satisfies $0 < v_P(J) < q+1$ --- so that the corresponding Drinfeld module acquires good reduction only over a ramified base --- we prove the bound
\[
v_P(a_{m,n})\;\geq\;\Big(v_P(J) + \tfrac1e\Big)\big(\psi(N) - m - n\big),
\; \mathrm{where}\;\; e = \tfrac{q+1}{\gcd(q+1,\,v_P(J))},
\]
for the coefficients of $\Phi_N(X+J, Y+J)$ (Theorem \ref{thm:ramified}); no CM hypothesis is required. The most interesting cases are nevertheless the singular moduli: we work out the modulus with CM by $\BF_q[\sqrt T\,]$ (the analogue of $D = -4$) and the modulus with CM by the conductor-$T$ order in $\BF_{q^2}[T]$ (the analogue of the non-maximal discriminants $D \in \{-12,-16,-27\}$ of \cite[Table 1]{Breuer2025}), for which the constants $\frac{q+2}2$ and $\frac{q+2}{q+1}$ are attained for every $N$ in our computations ($q=3$, resp.\ $q\in\{2,3\}$). As an application in the spirit of Wang's work on Kodaira types \cite{Wang2025}, $N$-isogenous $j$-invariants $J_1, J_2$ with $0 < v_P(J_1) < q+1$ at a prime $P\nmid N$ of degree one satisfy $v_P(J_2) = v_P(J_1)$, and even $v_P(J_2 - J_1) > v_P(J_1)$ when $J_1 \in A$ (Corollary \ref{cor:isogenyinvariance}).

\begin{rem}
Since $j=0$ is the only elliptic point of the Drinfeld modular line for $\mathrm{GL}_2(A)$, there is no analogue of the 
classical singular modulus $1728$: for every $J \neq 0$ the corresponding Drinfeld module has automorphism group $\BF_q^*$. 
\end{rem} 

In this study, we follow the strategy of \cite{Breuer2025}: an elementary interpolation lemma reduces the proof to obtaining lower bounds for $v_P\big(j(\tilde{\varphi}) - j(\tilde{\varphi}')\big)$ for suitable families of Drinfeld modules $\tilde{\varphi}, \tilde{\varphi}'$, and such bounds are provided by the theory of singular moduli. Here are some of the new ingredients:
\begin{itemize}
    \item The Gross--Zagier valuation formula \cite[Prop.~2.3]{GZ1985} is replaced by Dorman's exact formula \cite[Thm.~3.5]{Dorman1991} (see also Chen \cite{Chen2025} for arbitrary rank), for which we give a short self-contained proof, valid also in characteristic $2$ (Proposition \ref{GZ});
    \item the deformation theory of Serre--Tate and Lubin--Tate, used in \cite{Breuer2025} to construct congruence families of elliptic curves, is replaced by explicit manipulation of the coefficient $g$: every rank 2 Drinfeld module with good reduction has a model with $\Delta = 1$, and these models are parametrized by the affine $g$-line, a $(q+1)$-fold cover of the $j$-line ($j = g^{q+1}$) on which deformations can be written down by hand (Proposition \ref{LSTG});
    \item V\'elu's formulae are replaced by the elementary observation that a 
    cyclic $N$-isogeny with monic kernel polynomial $u = \sum_i u_i\tau^i$ transforms the coefficient $g$ by
    \[
    g' = u_0^{1-q}\,g + u_0^{-q}\,u_1\,(T^q - T),
    \]
    in which the factor $T^q - T = \prod_{\theta\in\BF_q}(T - \theta)$ --- the product of all monic primes of degree one --- leads to the estimate in Theorem \ref{thm:main}(1).
\end{itemize}

The paper is organized as follows. In Section \ref{sec:CJ} we introduce the counting functions $C_J(N,P)$ and describe them in the ordinary and supersingular cases. Section \ref{sec:proofs} contains the local theory and the proofs of the main theorems. Section \ref{sec:ramified} treats $j$-invariants which acquire good reduction only over a ramified extension, including two families of singular moduli. This section has no direct counterpart in \cite{Breuer2025}, where ramified bases enter only through explicit models of the thirteen rational singular moduli; here a single argument covers every $J$ with $0 < v_P(J) < q+1$, with or without complex multiplication. Section \ref{sec:computations} presents a new algorithm for computing $\Phi_N$ (Algorithm \ref{alg:phiN}). It runs the proof of Theorem \ref{thm:main} at finite precision: it deforms $\phi^0$ at a single place of $A$ and recovers the coefficients of $\Phi_N$ exactly from their $P$-adic expansions, without Chinese remaindering. Combined with the explicit height bound \eqref{explicitbound}, its output is provably correct (Proposition \ref{prop:algcorrect}). We then use it to verify our results numerically and to test the optimality of the constants. The code is available at \cite{GitHubRepo}.

\subsection{Use of AI}

The results in this paper are directly inspired by previous work of the authors. It is in particular in clear analogy with the recent \cite{Breuer2025} on elliptic modular polynomials, and it develops further the study of Drinfeld modular polynomials proposed in the series \cite{BPR2024, BRaz2024, BPP2026}. The current paper was started from
prompting Claude 5 Fable, with the aforementioned references as a thinking base, with GPT-5.5-Sol used to check correctness. Then Claude Opus 5.5 was used to revise explicit computations, and finally GPT-Astra was used to give feedback. The authors extensively rewrote the content to obtain the present text. The landscape evolves very quickly, they decided to carry this project as an experiment.

\section{The numbers $C_{J}(N,\Fp)$}\label{sec:CJ}

From now on, let $W$ be a complete discrete valuation ring, with valuation $v$, uniformizer $\pi$, fraction field $K$ of generic $A$-characteristic, and \emph{algebraically closed} residue field $\bar\BF$ of $A$-characteristic $\Fp = (P)$; we normalize $v(\pi) = 1$. The basic example is the completion $W = \hat{A^{\mathrm{ur}}_P}$ of the maximal unramified extension of $A_P$, with $\pi = P$.

Let $\varphi/K$ be a rank 2 Drinfeld module with good reduction, i.e.\ $\varphi$ is isomorphic over $K$ to a Drinfeld module $\varphi_T = \gamma(T) + g\tau + \Delta\tau^2$ with $g \in W$ and $\Delta \in W^*$; then its reduction $\bar \varphi$ modulo $\pi$ is again a rank 2 Drinfeld module. Let $J = j(\varphi)$ and $\CO_{J,\pi} = \End_{\bar\BF}(\bar \varphi)$. For monic $N$ with $P \nmid N$, define
\[
C_J(N,\pi) := \ord_X\big(\Phi_N(X + J, J) \bmod \pi\big),
\]
the number of cyclic kernels $C \subset \bar \varphi[N]$ with $j(\bar \varphi / C) = \bar J$, i.e.\ with $\bar \varphi/C\cong\bar \varphi$. More generally, for two such modules $\varphi^1, \varphi^2$ with $J_1 = j(\varphi^1)$ and $J_2 = j(\varphi^2)$, we put \[C_{J_1, J_2}(N,\pi) := \ord_X(\Phi_N(X+J_1, J_2)\bmod \pi),\]
the number of cyclic kernels $C\subset\bar \varphi^2[N]$ with $\bar \varphi^2/C\cong\bar \varphi^1$.

\begin{lemma}[Product formula]\label{lem:product}
    Let $N\in A$ be monic. The following properties hold.
    \begin{enumerate}
        \item[(i)] For every rank $2$ Drinfeld module $\varphi$ over an $A$-field $L$ of generic characteristic,
        \begin{equation}\label{eq:product}
            \Phi_N\big(X,\, j(\varphi)\big) \;=\; \prod_{C}\big(X - j(\varphi/C)\big) \qquad\text{in } \bar L[X],
        \end{equation}
        the product running over the $\psi(N)$ cyclic kernels $C\subset \varphi[N](\bar L)$.
        \item[(ii)] Suppose $\varphi$ is defined over $K$ with good reduction (a model with $g\in W$, $\Delta\in W^*$) and $v(N) = 0$. Then all $N$-torsion points of $\varphi$ lie in $W$, the nonzero ones being units with pairwise distinct reductions; every cyclic kernel $C\subset \varphi[N]$, its kernel polynomial $u_C$ and the quotient $\varphi/C$ are defined over $K$; each $\varphi/C$ has good reduction and $j(\varphi/C)\in W$; and $C\mapsto\bar C$ is a bijection onto the cyclic kernels of $\bar \varphi[N]$, with $\overline{\varphi/C} = \bar \varphi/\bar C$. Consequently, for any $J_1\in W$, reducing \eqref{eq:product} modulo $\pi$ shows that $C_{J_1,\,j(\varphi)}(N,\pi)$ is the number of cyclic kernels $\bar C\subset\bar \varphi[N]$ with $j(\bar \varphi/\bar C) = \bar J_1$.
    \end{enumerate}
\end{lemma}

\begin{proof}
    (1) Both sides are monic of degree $\psi(N)$ in $X$ with coefficients in the subfield $L_0\subseteq\bar L$ generated over $k$ by the coefficients of $\varphi$ and its (finitely many) $N$-torsion points. As $L_0$ has finite transcendence degree over $k$ while that of $\mathbb C_\infty$ is uncountable, there is a $k$-embedding $L_0\hookrightarrow\mathbb C_\infty$, and it suffices to prove \eqref{eq:product} over $\mathbb C_\infty$. There $\varphi\cong\phi^\Lambda$ for a rank $2$ $A$-lattice $\Lambda\subset\mathbb C_\infty$, where $\phi^\Lambda$ denotes the Drinfeld module attached to $\Lambda$, i.e.\ the one satisfying $e_\Lambda(az) = \phi^\Lambda_a\big(e_\Lambda(z)\big)$ for all $a\in A$, $e_\Lambda$ being the exponential function of $\Lambda$; thus $e_\Lambda$ induces an isomorphism of $A$-modules $\mathbb C_\infty/\Lambda\cong\phi^\Lambda(\mathbb C_\infty)$, and $\phi^\Lambda[N] \cong N^{-1}\Lambda/\Lambda$ \cite[Ch.~4]{Goss1996}. Under this isomorphism the cyclic kernels $C\subset \phi^\Lambda[N]$ correspond to the lattices $\Lambda\subseteq\Lambda_C\subseteq N^{-1}\Lambda$ with $\Lambda_C/\Lambda \cong A/NA$, and $\phi^\Lambda/C\cong\phi^{\Lambda_C}$, the isogeny being induced by the inclusion $\Lambda\subseteq\Lambda_C$. The identity \eqref{eq:product} is now the analytic construction of the Drinfeld modular polynomial: for $z$ on Drinfeld's upper half-plane and $\Lambda_z = Az + A$, the product $\prod_i\big(X - j(\Lambda_{z,i})\big)$ over the overlattices $\Lambda_{z,i}\supseteq\Lambda_z$ with $\Lambda_{z,i}/\Lambda_z\cong A/NA$ has coefficients that are modular functions, and equals $\Phi_N(X, j(z))$ for \emph{every} $z$; see \cite{Bae1992, Hsia1998}.

    (2) The division polynomial $\varphi_N(x)$ is $\BF_q$-linear with lowest coefficient $\gamma(N)\in W^*$ and leading coefficient a power of $\Delta\in W^*$; all its roots are integral over $W$, and its reduction is separable of the same degree $q^{2\deg N}$ (its derivative is $\overline{\gamma(N)}\neq0$). Since $K$ is complete and the residue field algebraically closed, Hensel's lemma lifts each root of the reduction uniquely to a root in $W$; comparing cardinalities, all $N$-torsion points of $\varphi$ lie in $W$, reduce bijectively onto $\bar \varphi[N]$, and the nonzero ones are units. The reduction map $\varphi[N]\to\bar \varphi[N]$ is additive and commutes with the $A$-actions, so it is an isomorphism of $A$-modules and induces a bijection on cyclic kernels. Kernel polynomials, being monic with coefficients the elementary symmetric functions of the points of $C$, have coefficients in $W$ and satisfy $\overline{u_C} = u_{\bar C}$; the quotient module $\varphi' := \varphi/C$, i.e.\ the Drinfeld module determined by $u_C\,\varphi_a = \varphi'_a\,u_C$ for all $a\in A$, then has coefficients in $W$ with $\Delta' = \Delta^{q^{\deg N}}\in W^*$ (compare coefficients in $u_C\,\varphi_T = \varphi'_T\,u_C$ from the top down), so $\varphi/C$ has good reduction, $j(\varphi/C)\in W$, and $\overline{\varphi/C} = \bar \varphi/\bar C$. The last claim follows by reducing \eqref{eq:product} with $X$ replaced by $X + J_1$.
\end{proof}

We shall also need dual isogenies.

\begin{lemma}[Dual isogenies]\label{lem:dual}
    Let $L$ be an algebraically closed $A$-field whose $A$-characteristic does not divide $N$, and let $u:\varphi\to\varphi'$ be a cyclic $N$-isogeny of rank $2$ Drinfeld modules over $L$ \textup{(}thus $u$ is separable, of $\tau$-degree $\deg N$\textup{)}. Then there is a unique $\hat u\in L\{\tau\}$ with $\hat u\,u = \varphi_N$. It is a cyclic $N$-isogeny $\hat u:\varphi'\to\varphi$ with $\ker\hat u = u(\varphi[N])\cong\varphi[N]/\ker u$, and $u\,\hat u = \varphi'_N$. In particular $\hat{\hat u} = u$, so $u\mapsto\hat u$ is a bijection from the set of cyclic $N$-isogenies $\varphi\to\varphi'$ onto the set of cyclic $N$-isogenies $\varphi'\to\varphi$.
\end{lemma}

\begin{proof}
    The additive polynomial $\varphi_N$ is separable, its derivative being $\gamma(N)\neq0$. Right division in $L\{\tau\}$ (possible since the leading coefficient of $u$ is a unit) gives $\varphi_N = \hat u\,u + r$ with $\deg_\tau r<\deg_\tau u = \deg N$. As $\ker u\cong A/NA$ is killed by $N$, the remainder $r$ vanishes on the $q^{\deg N}$ distinct roots of $u$, which forces $r = 0$. Next, $(\hat u\varphi'_a - \varphi_a\hat u)\,u = \hat u u\varphi_a - \varphi_a\varphi_N = 0$ for $a\in A$, so $\hat u$ is a morphism $\varphi'\to\varphi$ ($L\{\tau\}$ has no zero divisors), and likewise $(u\hat u - \varphi'_N)\,u = u\varphi_N - \varphi'_N u = 0$ gives $u\hat u = \varphi'_N$. Since $\hat u(u(x)) = \varphi_N(x) = 0$ for $x\in\varphi[N]$, we have $u(\varphi[N])\subseteq\ker\hat u$, and both sides have $|N|$ elements ($\hat u$ is separable, being a right factor of the separable $\varphi'_N$); hence $\ker\hat u = u(\varphi[N])\cong\varphi[N]/\ker u$. This quotient is cyclic: a generator of $\ker u\cong A/NA$ has, at each prime power $\ell^e\,\|\,N$, a coordinate which is a unit modulo $\ell$, so $\ker u$ is a direct summand of $\varphi[N]\cong(A/NA)^2$. Finally, $\hat{\hat u}$ is characterized by $\hat{\hat u}\,\hat u = \varphi'_N = u\,\hat u$, whence $\hat{\hat u} = u$.
\end{proof}

Recall \cite[Def.~4.1.9, Thm 4.1.3]{Papikian2023} that $\bar \varphi$ is \emph{supersingular} if $\End_{\bar\BF}(\bar \varphi)$ is an order in the quaternion algebra over $k$ ramified precisely at $P$ and $\infty$ (in which case it is a maximal order, by a theorem of Gekeler \cite{Gekeler1983}), and \emph{ordinary} otherwise (as we are in rank 2, there are no intermediate cases), in which case $\End_{\bar\BF}(\bar \varphi)$ is an $A$-order in an imaginary quadratic extension of $k$ (imaginary meaning that the place $\infty$ does not split).

\begin{prop}\label{prop:fullpsi}
    Suppose $P \nmid N$ and $\deg P = 1$. Then $C_0(N, P) = \psi(N)$.
\end{prop}

\begin{proof}
    The Drinfeld module $\phi^0$ with $\phi^0_T = \gamma(T) + \tau^2$ satisfies $j(\phi^0) = 0$ and has complex multiplication by $\BF_{q^2}[T]$, the integral closure of $A$ in the constant field extension $H = \BF_{q^2}(T)$ of $k$. A prime $P$ splits in $H/k$ if and only if $\deg P$ is even; hence for $\deg P$ odd the reduction of $\phi^0$ is supersingular \cite[Ex.~4.4.5 p.263]{Papikian2023}. For $\deg P = 1$, $j = 0$ is moreover the \emph{only} supersingular $j$-invariant in characteristic $P$ (see the formula in the odd case of \cite[Cor.~4.4.12]{Papikian2023}). All roots of $\Phi_N(X, 0)\bmod \pi$ are $j$-invariants of Drinfeld modules isogenous to the supersingular module $\bar\phi^0$, hence supersingular, hence equal to $0 \bmod \pi$. Since $\Phi_N(X,0)$ is monic of degree $\psi(N)$,
    the claim follows. 
\end{proof}

\begin{prop}\label{prop:deg2}
    Suppose $P\nmid N$ and $\deg P = 2$. There is exactly one supersingular $j$-invariant $j_P$ in characteristic $P$; it satisfies $j_P\in\BF_P := A/P$ and $j_P\neq0$. For every $J$ (the $j$-invariant of a Drinfeld module over $W$ with good reduction) with $J\equiv j_P\bmod\pi$ we have
    \[
    C_J(N,\pi) = \psi(N).
    \]
\end{prop}

\begin{proof}
    The number of supersingular $j$-invariants in characteristic $P$ equals $\frac{|P|-1}{q^2-1}$ when $\deg P$ is even and $\frac{|P|-q}{q^2-1}+1$ when $\deg P$ is odd (see \cite[\S5]{Gekeler1991}, \cite[\S4.4]{Papikian2023}), so there is exactly one for $\deg P = 2$. 
    Uniqueness forces $j_P$ to be fixed by $\mathrm{Gal}(\bar\BF_P/\BF_P)$, so $j_P\in\BF_P$; and $j_P \neq 0$ since $j = 0$ is ordinary in characteristic of even degree (by Proposition \ref{prop:fullpsi}'s proof). The equality $C_J(N,\pi) = \psi(N)$ follows exactly as in Proposition \ref{prop:fullpsi}: all modules isogenous to the (supersingular) reduction are supersingular, so all roots of $\Phi_N(X, J)\bmod\pi$ equal $j_P \equiv J$, and $P\nmid N$. 
\end{proof}

\begin{rem}
    Proposition \ref{prop:deg2} is the analogue of the classical fact that $C_5(N,13) = \psi(N)$, the invariant $5\in\BF_{13}$ being the unique supersingular $j$-invariant in characteristic $13$ \cite[Prop.~2.1]{Breuer2025}. A concrete example with $q = 3$: the unique supersingular invariant in characteristic $P = T^2+T+2$ is $j_P \equiv T + 2 \bmod P$, the reduction of the singular modulus $J_s = T^2(1+T)^4$ of Section \ref{sec:ramified}. Note that, in contrast with the classical case, $j_P$ need not lie in the prime field $\BF_p$.
\end{rem}

\begin{prop}\label{ordinary}
    Suppose $P\nmid N$ and $\varphi$ has ordinary reduction; then $\CO_{J,\pi}=\End_{\bar\BF}(\bar \varphi)$ is an $A$-order in an imaginary quadratic extension of $k$.
    \begin{enumerate}
        \item We have
        \[
        C_J(N,\pi) \leq \#I_{\mathrm{cyc}}(N, \CO_{J,\pi}) := \#\{\mathfrak{n}\subset\CO_{J,\pi}\;\text{ideal} : \CO_{J,\pi}/\mathfrak{n} \cong A/NA\},
        \]
        with equality if $\CO_{J,\pi}$ is a principal ideal domain.
        \item If $\varphi/K$ has complex multiplication and $P$ does not divide the conductor of $R=\End_{\bar{K}}(\varphi)$, then
        \[
        C_J(N,\pi) = C_J(N, 0) := \ord_X\big(\Phi_N(X+J, J) \in K[X]\big).
        \]
    \end{enumerate}
\end{prop}

\begin{proof}
    We prove (1) using the theory of ordinary Drinfeld modules over finite fields: a kernel $C$ with $\bar \varphi/C\cong\bar \varphi$ is the kernel $\bar \varphi[(f)]$ of an endomorphism $f\in\CO=\CO_{J,\pi}$ of norm in $\BF_q^* N$, so $C = \bar \varphi[\mathfrak n]$ for a \emph{principal} ideal $\mathfrak n = (f) \in I_{\mathrm{cyc}}(N,\CO)$, and conversely; distinct such ideals have distinct kernels. If $\CO$ is a principal ideal domain, all ideals in $I_{\mathrm{cyc}}(N,\CO)$ arise. (Alternatively, (1) follows from the simply transitive action of $\mathrm{Pic}(\CO)$ on the isomorphism classes with endomorphism ring $\CO$, see \cite[Cor.~p.170]{Yu1995}.)

    For part (2), suppose $\varphi$ has CM by the order $R = \End_{\bar K}(\varphi)$ of conductor $\mathfrak c$, and let $C\subset \varphi[N]$ be a cyclic kernel, defined over $K$, with $j(\varphi/C) \equiv J \bmod \pi$; we must show $j(\varphi/C) = J$. Both $\varphi/C$ and $\varphi$ have CM by orders in the same field, and since $P\nmid N\mathfrak{c}$, the conductor of $\End(\varphi/C)$ is not divisible by $P$. Its prime-to-$P$ part agrees with that of $\mathfrak c$ as well, because both equal the prime-to-$P$ part of the conductor of the common reduced endomorphism ring $\CO_{J,\pi}$ (see \cite[Prop.~2.1]{BPP2026}). Hence $\End(\varphi/C) \cong R$, and now $j(\varphi/C)\equiv J\bmod\pi$ forces $j(\varphi/C) = J$ because reduction is injective on the $j$-invariants of CM Drinfeld modules with a fixed ordinary endomorphism ring — this is part of Theorem D of \cite{BPP2026} (due to Bae and Koo), where injectivity comes from the simply transitive action of the class group \cite{Yu1995} and surjectivity from the Deuring lifting theorem for Drinfeld modules \cite[Thm.~7]{CP2015}. Thus every kernel counted by $C_J(N,\pi)$ is counted by $C_J(N,0)$ (and kernels of $\varphi[N]$ and $\bar \varphi[N]$ correspond bijectively since $P\nmid N$); the converse inclusion is clear, whence equality.
\end{proof}

\begin{rem}
    In case (2), the coefficient $a_{C_J(N,\pi),0}$ of $\Phi_N(X+J,J)$ is a $\pi$-adic unit, so no divisibility statement of the shape of Theorem \ref{thm:main} can hold at ordinary primes not dividing the conductor. Theorem \ref{thm:main} only gives non-trivial divisibility of coefficients when $C_J(N,\pi) > C_J(N,0)$. Indeed, we shall prove in Remark \ref{rem:vanishing} that the coefficients $a_{m,n}$ of $\Phi_N(X+J,Y+J)$ with $m+n < C_J(N,0)$ all vanish.
\end{rem}

\begin{prop}\label{supersingular}
    Let $P\nmid N$ and suppose $\varphi$ has supersingular reduction, so that $\CO = \CO_{J,\pi}$ is a maximal order in the quaternion algebra $D$ over $k$ ramified exactly at $P$ and $\infty$. The following properties hold.
    \begin{enumerate}
        \item We have
        \[
        C_J(N,\pi) = \frac{1}{\#\CO^*}\sum_{\substack{d \text{ monic}\\ d^2 \mid N}}\mu(d)\,
        \#\big\{f\in\CO \;:\; \nrd(f) \in \BF_q^*\cdot N/d^2\big\},
        \]
        where $\mu$ denotes the M\"obius function for monic polynomials and $\nrd$ the reduced norm on $D$.
        \item Suppose $q$ is odd, and suppose moreover that $\varphi/K$ has complex multiplication by an imaginary quadratic $A$-order of discriminant $\mathfrak d$, with $P\nmid\mathfrak d$, and that $C_J(N,\pi) > C_J(N,0)$. Then
        \[
        \deg P \;\leq\; \deg \mathfrak{d} + \deg N.
        \]
    \end{enumerate}
\end{prop}

For the proof of part (2) we need the following analogue of a lemma of Kaneko \cite{Kaneko1989}.

\begin{lemma}\label{lem:kaneko}
    Let $q$ be odd and let $\CO$ be a maximal order in the quaternion algebra $D/k$ ramified exactly at $P$ and $\infty$. Suppose $s_1, s_2\in\CO$  with $\mathrm{trd}(s_i) = 0$, $s_i^2 = d_i\in A$, which do not commute. Then $t := \mathrm{trd}(s_1s_2)\in A$ satisfies
    \[
    0 \neq t^2 - 4d_1d_2 \equiv 0 \bmod P, \qquad \deg(t^2 - 4d_1d_2)\leq \deg d_1 + \deg d_2;
    \]
    in particular $\deg P \leq \deg d_1 + \deg d_2$.
\end{lemma}

\begin{proof}
    Set $w := s_1s_2 - s_2s_1 \neq 0$. Using $s_2s_1 = t - s_1s_2$ one computes $(s_1s_2)^2 = ts_1s_2 - d_1d_2$ and hence
    $w^2 = (2s_1s_2 - t)^2 = t^2 - 4d_1d_2 \in A\setminus\{0\}$.
    The $A$-lattice $\Lambda = A + As_1 + As_2 + As_1s_2 \subseteq \CO$ has rank $4$ (as $1, s_1, s_2, s_1s_2$ are linearly independent when $s_1, s_2$ do not commute), and the Gram matrix of the reduced trace form in this basis is
    \[
    \begin{pmatrix}
        2 & 0 & 0 & t\\ 0 & 2d_1 & t & 0\\ 0 & t & 2d_2 & 0 \\ t & 0 & 0 & t^2 - 2d_1d_2
    \end{pmatrix},
    \]
    with determinant $-(t^2 - 4d_1d_2)^2$. Thus the reduced discriminant of $\Lambda$ is generated by $t^2 - 4d_1d_2$, while the reduced discriminant of the maximal order $\CO$ is $(P)$ (see \cite[\S41]{Voight2021} for the general theory, which applies verbatim over $A$; the finite ramified prime of $D$ is $P$). Since $\Lambda\subseteq\CO$, the discriminant of $\CO$ divides that of $\Lambda$, so we get $P \mid t^2 - 4d_1d_2$.

    For the degree bound, note that $k(w)$ is a commutative subalgebra of the division algebra $D\otimes_k k_\infty$ (as $D$ is ramified at $\infty$), hence a field; so $w^2 = t^2 - 4d_1d_2$ is a non-square in $k_\infty = \BF_q((1/T))$. Since squares in $k_\infty^*$ are exactly the elements of even degree whose leading coefficient is a square in $\BF_q^*$, and since $\deg (t^2) > \deg (4d_1d_2)$ would force $t^2 - 4d_1d_2$ to have even degree equal to $2\deg t$ and square leading coefficient (which is excluded by the previous sentence), we conclude $\deg (t^2) \leq \deg(d_1d_2)$ and therefore $\deg(t^2 - 4d_1d_2) \leq \deg d_1 + \deg d_2$.
\end{proof}

We will also make use of the following lemma.
\begin{lemma}\label{lem:CMtate}
    Let $\varphi/K$ be a rank $2$ Drinfeld module of generic $A$-characteristic whose endomorphism ring $R := \End_{\bar K}(\varphi)$ is an $A$-order in an imaginary quadratic field $F$. Let $\ell$ be a monic prime of $A$ and write $R_\ell := R\otimes_AA_\ell$. Then the Tate module $T_\ell \varphi$ is isomorphic, as an $R_\ell$-module, to an $R_\ell$-lattice in $F_\ell := F\otimes_kk_\ell$ whose multiplicator ring is exactly $R_\ell$. Consequently
    \[
    \End_{R_\ell}(T_\ell \varphi) \;=\; R_\ell,
    \]
    i.e.\ every $R_\ell$-linear endomorphism of $T_\ell \varphi$ is multiplication by an element of $R_\ell$.
\end{lemma}

\begin{proof}
    We may embed a field of definition of $\varphi$ and its endomorphisms into $\mathbb C_\infty$; this changes neither the torsion $A$-modules $\varphi[\ell^m]$ nor $\End$ (an endomorphism is an $\BF_q$-linear polynomial determined by its values on the infinitely many torsion points, and its coefficients are algebraic over the field of definition). Over $\mathbb C_\infty$ we have $\varphi\cong\phi^\Lambda$ for a rank $2$ $A$-lattice $\Lambda$ (notation as in the proof of Lemma \ref{lem:product}); then $\End(\varphi)\cong\{x\in\mathbb C_\infty : x\Lambda\subseteq\Lambda\}$ and the exponential function $e_\Lambda$ identifies $\varphi[\ell^m]$ with $\ell^{-m}\Lambda/\Lambda$ \cite[Ch.~4]{Goss1996}.
    Thus $\Lambda$ is a module over its multiplicator ring $R$, torsion-free of rank one (as $\Lambda\otimes_Ak$ is a $1$-dimensional $F$-vector space), hence isomorphic to a fractional $R$-ideal of $F$ with multiplicator ring exactly $R$, and
    \[
    T_\ell \varphi \;=\; \varprojlim_m\,\ell^{-m}\Lambda/\Lambda \;\cong\; \Lambda\otimes_AA_\ell.
    \]
    Multiplicator rings of lattices commute with the flat base change $A\to A_\ell$, so $T_\ell \varphi$ is an $R_\ell$-lattice in $F_\ell$ with multiplicator ring $R_\ell$. Finally, an $R_\ell$-linear endomorphism of $T_\ell \varphi$ extends $k_\ell$-linearly to $T_\ell \varphi\otimes_{A_\ell}k_\ell\cong F_\ell$, where it is $F_\ell$-linear, hence multiplication by an element of $F_\ell$; as it stabilizes $T_\ell \varphi$, that element lies in the multiplicator ring $R_\ell$.
\end{proof}

\begin{proof}[Proof of Proposition \ref{supersingular}] Let us start with part (1).
    Cyclic kernels $C\subset \bar \varphi[N]$ with $\bar \varphi/C \cong \bar \varphi$ correspond to endomorphisms $f \in \CO$ with cyclic kernel and $\nrd(f) \in \BF_q^*N$, each kernel arising from exactly $\#\CO^*$ endomorphisms (for each post-composition with an automorphism). Thus
    \[
    \#\{f \in \CO : \nrd(f)\in\BF_q^* N, \; \ker f \text{ cyclic}\} = \#\CO^* \cdot C_J(N,\pi).
    \]
    An arbitrary $f\neq0$ with $\nrd(f) \in \BF_q^*N$ factors uniquely as $f = f'\circ \varphi_d$, 
    where $d$ is monic with $d^2 \mid N$ and $f'$ has cyclic kernel; M\"obius inversion gives the formula. (One has $\#\CO^* \in \{q-1, q^2-1\}$.)

    Let us now prove (2). The first non-vanishing coefficient of $\Phi_N(X + J, J) \in K[X]$ is
    $a_{C_J(N,0),0} = \pm\prod \big(j(\varphi/C) - J\big)$,
    the product over cyclic kernels $C$ with $j(\varphi/C) \neq J$.
    If $C_J(N,\pi) > C_J(N,0)$, this product has positive valuation, so some cyclic kernel $C$ has $j(\varphi/C)\neq J$ but $\overline{\varphi/C} \cong \bar \varphi$. Composing $\varphi\to \varphi/C$ with such an isomorphism yields $ f\in\CO$ with $\nrd(f)\in\BF_q^*N$ and $\ker f = \bar C$. Let $s_1\in\CO$ be the image under reduction of a trace-zero generator of the CM order, so $s_1^2 = d_1$ with $(d_1) = (\mathfrak d)$ up to squares of units and $\deg d_1 = \deg\mathfrak d$; reduction of endomorphisms is injective, and $s_1\notin A$. (Injectivity: let $u = \sum_{i\leq n}u_i\tau^i$, $u_n\neq0$, be a non-zero endomorphism of $\varphi$ over $K$. Comparing the coefficients of $\tau^{n+2}$ in $u\varphi_T = \varphi_Tu$ gives $u_n^{q^2-1} = \Delta^{q^n-1}\in W^*$, so $u_n\in W^*$; and descending induction gives $u_i\in W$ for all $i$, because the coefficients of $\tau^{i+2}$ give $\Delta u_i^{q^2} - \Delta^{q^i}u_i\in W$ once $u_{i+1}, u_{i+2}\in W$, and $\Delta X^{q^2} - aX\in W$ with $a,\Delta\in W^*$ forces $X\in W$. Thus $u$ reduces to an endomorphism $\bar u$ of $\bar\varphi$ whose leading coefficient $\bar u_n$ is non-zero, so $\bar u\neq0$.)

    Suppose first that $f$ commutes with $s_1$, and write $R = \End_K(\varphi)$, an order of discriminant $\mathfrak d$, with reduction $\bar R\subseteq\CO$; note that $R = \End_{\bar K}(\varphi)$, since an endomorphism is determined by its values on the prime-to-$P$ torsion, which is $K$-rational by Lemma \ref{lem:product}(2). Then $f$ lies in the $A$-order $S := \CO\cap k(s_1) \supseteq \bar R$. We claim $S = \bar R$. At the prime $P$: since $P\nmid\mathfrak d$, the order $\bar R\otimes A_P$ is maximal, so $S\otimes A_P = \bar R\otimes A_P$. At a prime $\ell\neq P$: reduction is a bijection on $\ell$-power torsion (Lemma \ref{lem:product}(2), applied with $N = \ell^m$), inducing an isomorphism $T_\ell \varphi \xrightarrow{\;\sim\;} T_\ell\bar \varphi$ of Tate modules compatible with reduction of endomorphisms. Every element $x\in S$ commutes with $\bar R$ and stabilizes $T_\ell\bar \varphi$; transporting through this isomorphism, $x$ acts on $T_\ell \varphi$ as an $R\otimes A_\ell$-linear endomorphism, hence, by Lemma \ref{lem:CMtate}, as multiplication by some $r\in R\otimes A_\ell$. Since $\CO\otimes A_\ell$ acts faithfully on $T_\ell\bar \varphi\otimes_{A_\ell}k_\ell$ (the quaternion algebra $D\otimes_kk_\ell$ is simple and acts nontrivially), $x$ equals the image $\bar r\in\bar R\otimes A_\ell$; thus $S\otimes A_\ell \subseteq \bar R\otimes A_\ell$, and therefore $S\otimes A_\ell = \bar R\otimes A_\ell$ for all $\ell\neq P$. Whence $S = \bar R$ and $f = \bar f_0$ for a unique $f_0\in R$. Since $P\nmid N$, torsion reduces isomorphically and $C = \ker(f_0)\cap \varphi[N]$, so that $\varphi/C\cong f_0(\varphi) = \varphi$ and $j(\varphi/C) = J$, contradicting $j(\varphi/C)\neq J$. Therefore some kernel yields $f$ \emph{not} commuting with $s_1$.

    For such $f$, put $s_2 := f - \tfrac{1}{2}\mathrm{trd}(f)$, so $\mathrm{trd}(s_2) = 0$ and $s_2^2 = d_2 := \tfrac14(\mathrm{trd}(f)^2 - 4\nrd(f))$. The subfield $k(f)\subset D$ is imaginary (it embeds in the division algebra $D\otimes k_\infty$), so as in the proof of Lemma \ref{lem:kaneko}, $\deg \mathrm{trd}(f)^2 \leq \deg 4\nrd(f) = \deg N$, whence $\deg d_2 \leq \deg N$. Lemma \ref{lem:kaneko} applied to $s_1, s_2$ gives $\deg P \leq \deg d_1 + \deg d_2 \leq \deg\mathfrak d + \deg N$.
\end{proof}

\section{Proof of the main results}\label{sec:proofs}

\subsection{Interpolation}

We shall use the following elementary lemma, which is \cite[Lemma 3.3]{Breuer2025}; its proof (via Vandermonde determinants) is independent of the ground field and we do not repeat it.

\begin{lemma}\label{lem:interpolation}
    Let $f(Y) = a_0 + a_1Y + \cdots + a_dY^d\in W[Y]$. Fix $n \in\BZ$ and let $y_0,\ldots,y_d\in W$ satisfy
    $v(y_m) = n$ for all $m$ and $v(y_m - y_l) = n$ for all $m\neq l$. Then for any $s\in\{ 0, 1, \ldots, d\}$,
    \[
    v(a_s) \geq \min_{0\leq m\leq d} v(f(y_m)) - ns.
    \]
\end{lemma}

\begin{prop}\label{prop:general}
    Let $\varphi^1, \varphi^2$ be rank 2 Drinfeld modules over $W$ with good reduction and $j$-invariants $J_1, J_2 \in W$. Let $N$ be monic with $v(N) = 0$. Fix positive integers $n_1, n_2$ and suppose there exist rank 2 Drinfeld modules $\tilde{\varphi}^0, \tilde{\varphi}^1, \ldots, \tilde{\varphi}^{\psi(N)}$ over $W$ such that:
    \begin{enumerate}
        \item for each $s\in\{0,\ldots, \psi(N)\}$, $\tilde{\varphi}^s$ has good reduction;
        \item $v(j(\tilde{\varphi}^s) - J_2) = v(j(\tilde{\varphi}^s) - j(\tilde{\varphi}^l)) = n_2$ for all $s \neq l$;
        \item for every $s$ and every $\tilde{\varphi}'^s$ linked to $\tilde{\varphi}^s$ by a cyclic $N$-isogeny,
        \[
        v(j(\tilde{\varphi}'^s) - J_1) > 0 \implies v(j(\tilde{\varphi}'^s) - J_1)\geq n_1.
        \]
    \end{enumerate}
    Then the coefficients of $\Phi_N(X + J_1, Y + J_2) = \sum_{m,n}a_{m,n}X^mY^n \in W[X,Y]$ satisfy
    \[
    v(a_{m,n}) \;\geq\; n_1 C_{J_1,J_2}(N,\pi) - n_1 m - n_2 n.
    \]
\end{prop}

\begin{proof}
    We follow the proof of \cite[Prop.~3.4]{Breuer2025}.
    Write 
    \[
    \Phi_N(X + J_1, Y+J_2) = \sum_{m=0}^{\psi(N)} b_m(Y)X^m, \qquad b_m(Y) = \sum_n a_{m,n}Y^n.
    \]
    Fix $s$ and let $y_s = j(\tilde{\varphi}^s) - J_2$. Since $v(N) = 0$ and $\tilde{\varphi}^s$ has good reduction, Lemma \ref{lem:product}(2) applies to $\tilde{\varphi}^s$: all $N$-torsion points, all cyclic kernels $C\subset\tilde{\varphi}^s[N]$ and all quotient modules $\tilde{\varphi}^s/C$ are defined over $W$ and have good reduction.

    By Lemma \ref{lem:product}, $\Phi_N(X, y_s + J_2) = \prod_C\big(X - j(\tilde{\varphi}^s/C)\big)$. Since $y_s + J_2 \equiv J_2\bmod\pi$, we have $\Phi_N(X + J_1, y_s + J_2)\equiv\Phi_N(X+J_1,J_2)\bmod\pi$, so exactly $C_{J_1,J_2}(N,\pi)$ of the numbers $j(\tilde{\varphi}^s/C) - J_1$ (counted with multiplicity) have positive valuation; by hypothesis (3) these have valuation $\geq n_1$, while the remaining $\psi(N) - C_{J_1,J_2}(N,\pi)$ are units. Hence $b_m(y_s)$, being up to sign the $(\psi(N)-m)$-th elementary symmetric polynomial in the $j(\tilde{\varphi}^s/C) - J_1$, satisfies
    \[
    v(b_m(y_s)) \geq n_1(C_{J_1,J_2}(N,\pi) - m), \qquad 0\leq m \leq C_{J_1,J_2}(N,\pi).
    \]
    Applying Lemma \ref{lem:interpolation} to each $b_m$ with the interpolation points $y_s$ (hypothesis (2)) gives $v(a_{m,n}) \geq n_1(C_{J_1,J_2}(N,\pi)-m) - n_2 n$ for $m\leq C_{J_1,J_2}(N,\pi)$, and for $m > C_{J_1,J_2}(N,\pi)$ the claimed bound is vacuous or follows from $v(a_{m,n})\geq 0$.
\end{proof}

\subsection{Valuations of differences of $j$-invariants}

The following is the rank 2 Drinfeld module analogue of the Gross--Zagier valuation formula \cite[Prop.~2.3]{GZ1985}. It is due to Dorman \cite[Thm.~3.5]{Dorman1991}, who works in odd characteristic throughout. We include the short proof, which works in every characteristic (this matters for our computations with $q = 2$), both for completeness and because we shall use the auxiliary statements it contains. Closely related results on singular and supersingular moduli of Drinfeld modules appear in Brown \cite{Brown1992}, and more recently Chen \cite[Prop.~3.11]{Chen2025} proves an inequality in arbitrary rank. 

\begin{prop}\label{GZ}
    Let $\varphi, \varphi'$ be rank 2 Drinfeld modules over $W$ with good reduction and $j = j(\varphi) \neq j' = j(\varphi')$. Then
    \[
    v(j - j') = \frac{1}{q-1}\sum_{m\geq 1}\#\Isom_{W/\pi^m}(\varphi, \varphi'),
    \]
    where $\Isom_{W/\pi^m}$ denotes the set of isomorphisms of Drinfeld modules over $W/\pi^m$ between the reductions.
\end{prop}

\begin{proof}
    Since the residue field $\bar\BF$ is algebraically closed and $\Delta, \Delta'\in W^*$ are units, Hensel's lemma allows us to solve $c^{q^2-1} = \Delta^{-1}$ in $W^*$ and replace $\varphi$ by the isomorphic module $c^{-1}\varphi c$; thus we may and do assume
    \[
    \varphi_T = \gamma(T) + g\tau + \tau^2, \qquad \varphi'_T = \gamma(T) + g'\tau + \tau^2, \qquad g, g' \in W,
    \]
    so that $j = g^{q+1}$ and $j' = g'^{q+1}$. As $q^2 - 1$ is prime to the residue characteristic, the group of $(q^2-1)$th roots of unity $\mu_{q^2-1}(W/\pi^m)$ maps isomorphically onto $\mu_{q^2-1}(\bar\BF)\cong\BF_{q^2}^*$ for every $m\geq1$. An isomorphism $\varphi\to\varphi'$ over $W/\pi^m$ is an element $c \in (W/\pi^m)^*$ with $c^{q^2-1} = 1$ and $c^{1-q}g \equiv g' \bmod \pi^m$; writing $\zeta = c^{1-q}$, which runs over $\mu_{q+1}$ with each value attained by exactly $q-1$ elements $c \in \mu_{q^2-1}$, we obtain
    \begin{equation}\label{isommu}
    \#\Isom_{W/\pi^m}(\varphi,\varphi') = (q-1)\cdot\#\{\zeta\in\mu_{q+1} : v(g' - \zeta g)\geq m\}.
     \end{equation}
    Summing over $m \geq 1$ and using $x^{q+1} - y^{q+1} = \prod_{\zeta\in\mu_{q+1}}(x - \zeta y)$, we conclude
    \[
    \frac{1}{q-1}\sum_{m\geq1}\#\Isom_{W/\pi^m}(\varphi,\varphi') = \sum_{\zeta\in\mu_{q+1}}v(g' - \zeta g) = v\big(g'^{q+1} - g^{q+1}\big) = v(j - j').\qedhere
    \]
\end{proof}

\begin{corollary}\label{cor:lattice}
    Let $\varphi$ be a rank 2 Drinfeld module over $W$ with good reduction and $j(\varphi)\neq 0$. Then
    \[
    v(j(\varphi)) \in (q+1)\,\BZ_{\geq 0}.
    \]
    Conversely, if $J \in W$ satisfies $(q+1) \mid v(J)$, then there is a Drinfeld module over $W$ with good reduction and $j$-invariant $J$.
\end{corollary}

\begin{proof}
    With the normalization of the previous proof, if $j(\varphi) = g^{q+1}$ with $g \in W$, we have $v(j) = (q+1)v(g)$. Conversely, take $g$ to be any element satisfying $g^{q+1} = J$; if $(q+1)\mid v(J)$ then, after factoring out $\pi^{v(J)/(q+1)}$, Hensel's lemma provides such a $g$ in $W$, and $\varphi_T = \gamma(T) + g\tau + \tau^2$ has good reduction.
\end{proof}

Corollary \ref{cor:lattice} is the exact analogue of the classical facts that $v(j) \equiv 0 \bmod 3$ and $v(j - 1728)\equiv 0\bmod 2$ for elliptic curves with good reduction over an absolutely unramified base --- except that in the Drinfeld setting a single congruence controls everything, and there is no wild deepening at small primes. We also record the automorphism counts contained in the proof of Proposition \ref{GZ} (see also \cite[Prop.~3.1]{Dorman1991} in odd characteristic), which we will use in the next subsection.

\begin{corollary}\label{cor:autcount}
    Let $\varphi_T = \gamma(T) + g\tau + \tau^2$ over $W$. Then for $m \geq 1$:
    \[
    \#\Aut_{W/\pi^m}(\varphi) = \begin{cases}
        q^2 - 1 & \text{if } v(g)\geq m,\\
        q - 1 & \text{otherwise.}
    \end{cases}
    \]
\end{corollary}
\begin{proof}
    This follows from (\ref{isommu}).
\end{proof}

\subsection{Deformations of Drinfeld modules}

In this paragraph we record a result on deformations of a given Drinfeld module. 

\begin{prop}\label{LSTG}
    Let $\varphi/W$ be a rank 2 Drinfeld module with good reduction and let $M\geq1$. Then there exist infinitely many Drinfeld modules $\tilde{\varphi}^1, \tilde{\varphi}^2, \ldots$ over $W$ such that
    \begin{enumerate}
        \item[(i)] for any integer $s\ge 1$ we have $\tilde{\varphi}^s \equiv \varphi \bmod \pi^M$ (in particular each $\tilde{\varphi}^s$ has good reduction), and
        \item[(ii)] for all integers $s,l$ with $s\neq l$,
        \begin{equation}\label{valautom}
        v(j(\tilde{\varphi}^s) - j(\varphi)) = v(j(\tilde{\varphi}^s) - j(\tilde{\varphi}^l)) = \frac{1}{q-1}\sum_{m=1}^{M}\#\Aut_{W/\pi^m}(\varphi).
         \end{equation}
    \end{enumerate}
\end{prop}

\begin{proof}
    Normalize $\varphi : \varphi_T = \gamma(T) + g\tau + \tau^2$ as in Proposition \ref{GZ} and set
    \[
    \tilde{\varphi}^s : \; T \mapsto \gamma(T) + (g + \eps_s\pi^M)\tau + \tau^2,
    \]
    where the units $\eps_s\in W^*$ are chosen as follows. By Corollary \ref{cor:autcount}, the right-hand side of (\ref{valautom}) equals $M(q+1)$ if $v(g)\geq M$, and $v(g)(q+1) + (M - v(g))$ if $v(g) < M$. All these modules have $\Delta = 1$, so their $j$-invariants are the $(q+1)$-st powers of their coefficients of $\tau$. The factorization $x^{q+1} - y^{q+1} = \prod_{\zeta\in\mu_{q+1}}(x - \zeta y)$ used in the proof of Proposition \ref{GZ} (recall that $\mu_{q+1}\subset W^*$) therefore gives
    \[
    v(j(\tilde{\varphi}^s) - j(\tilde{\varphi}^l)) = \sum_{\zeta\in\mu_{q+1}}v\big((g + \eps_s\pi^M) - \zeta(g + \eps_l\pi^M)\big),
    \]
    and similarly for $v(j(\tilde{\varphi}^s) - j(\varphi))$ with $\eps_l = 0$. We compute the valuation of each factor $(1-\zeta)g + (\eps_s - \zeta\eps_l)\pi^M$.
    For $\zeta = 1$ the factor is $(\eps_s - \eps_l)\pi^M$, of valuation $M$ provided $\bar\eps_s\neq\bar\eps_l$.
    For $\zeta\neq1$ (so that $1 - \zeta$ is a unit): if $v(g) < M$ the factor has valuation $v(g)$ unconditionally; if $v(g) > M$ it has valuation $M$ provided $\bar\eps_s \neq \zeta\bar\eps_l$; and if $v(g) = M$ it has valuation $M$ provided $\bar\eps_s - \zeta\bar\eps_l + \overline{(1-\zeta)g\pi^{-M}}\neq 0$.
    These constitute finitely many polynomial conditions on the pair $(\bar\eps_s, \bar\eps_l)\in(\bar\BF^*)^2$, none identically satisfied; since $\bar\BF$ is infinite we may choose a sequence $(\bar\eps_s)_s$ avoiding all of them. Summing the factor valuations gives $M + q\min(v(g), M)$, which equals the right-hand side of (\ref{valautom}) by Corollary \ref{cor:autcount}.
\end{proof}

\subsection{Isogenies and the coefficient $g$}

Let $\varphi_T = \gamma(T) + g\tau + \Delta\tau^2$ over $W$ with good reduction, let $C \subset \varphi[N]$ be a cyclic kernel (with $v(N) = 0$) and let
\[
u_C(x) = \prod_{c\in C}(x - c) = \sum_{i=0}^{d} u_i\tau^i, \qquad d = \deg N, \quad u_d = 1,
\]
be the associated (monic, linearized) kernel polynomial, so that $u_C$ is a cyclic $N$-isogeny $\varphi\to\varphi' := \varphi/C$. Here the quotient $\varphi/C$ is the Drinfeld module determined by $u_C\,\varphi_a = \varphi'_a\,u_C$ for all $a\in A$; it exists because $u_C\varphi_a$ vanishes on $C$ and is therefore right divisible by $u_C$ (as in the proof of Lemma \ref{lem:dual}). Every cyclic $N$-isogeny from $\varphi$ arises this way up to isomorphism. Comparing coefficients of $\tau^{d+2}$, $\tau^{d+1}$ and $\tau^{1}$ in $u_C\varphi_T = \varphi'_T u_C$ yields
\begin{align}
    \Delta' &= \Delta^{q^d},\nonumber\\
    g' &= g^{q^d} + u_{d-1}\Delta^{q^{d-1}} - \Delta^{q^d}u_{d-1}^{q^2},\label{eq:gprime-top}\\
    g' &= u_0^{1-q}\, g + u_0^{-q}\, u_1\,\big(\gamma(T)^q - \gamma(T)\big).\label{eq:gprime-bottom}
\end{align}
Both expressions for $g'$ will be used: \eqref{eq:gprime-top} shows integrality, while \eqref{eq:gprime-bottom} is an exact formula isolating the arithmetic. Note that all torsion points $c\in C\setminus\{0\}$ are units by Lemma \ref{lem:product}(ii). Moreover, if we denote by 
\begin{equation}\label{SC}
    S_C := \sum_{0\neq c\in C}c^{1-q},
\end{equation}
    we have $u_1/u_0 = S_C$ by comparing the expansions of $u_C'(x)/u_C(x) = \sum_{c}(x-c)^{-1}$, so
\[
v(u_0) = 0, \qquad u_1 = u_0\, S_C.
\]
 Finally, in $A$ we have the factorization
\[
T^q - T = \prod_{\theta\in\BF_q}(T - \theta),
\]
the product of \emph{all monic primes of degree one}, each to the first power. Hence
\begin{equation}\label{eq:vTqT}
    v\big(\gamma(T)^q - \gamma(T)\big) = \begin{cases} e & \text{if } \deg P = 1,\\ 0 & \text{if } \deg P \geq 2,\end{cases}
\end{equation}
where $e = v(P)$ is the ramification index of $W$ over $A_P$.

At primes of degree one, \eqref{eq:gprime-bottom} has the following consequence for the kernel polynomials of modules with supersingular reduction $j = 0$, which will be the key to the fractional exponents of Section \ref{sec:ramified}.

\begin{lemma}\label{lem:Fq2}
    Let $P = T-\theta$ be a prime of degree one, $P\nmid N$, and let $\varphi$ be a rank 2
    Drinfeld module over $W$ whose reduction is $\bar\varphi_T = \theta + \tau^2$.
    Then the reduction $\bar u_C$ of the kernel polynomial of any cyclic kernel
    $C\subset\varphi[N]$ has coefficients in $\BF_{q^2}$; in particular
    $\bar u_0\in\BF_{q^2}^*$.
\end{lemma}
\begin{proof}
    $\bar u_C$ is the kernel polynomial of $\bar C\subset\bar\varphi[N]$, and by
    \eqref{eq:gprime-bottom} the quotient module $\bar\varphi/\bar C$
    has coefficients $\bar g' = \bar u_0^{-q}\bar u_1(\theta^q - \theta) = 0$
    (as $\theta\in\BF_q$) and $\bar\Delta' = 1$; hence
    $\bar u_C(\theta + \tau^2) = (\theta+\tau^2)\bar u_C$, i.e.\ $\bar u_C$
    commutes with $\tau^2$. Writing out $\tau^2 a = a^{q^2}\tau^2$, this forces
    all coefficients of $\bar u_C$ to lie in $\BF_{q^2}$. ($\bar u_0 \neq 0$
    because the torsion points are units, Lemma \ref{lem:product}(ii).)
\end{proof}

\subsection{Proof of Theorem \ref{thm:main}}

We first prove the general local statement.

\begin{thm}\label{thm:bigprimes}
    Let $\varphi^1, \varphi^2$ be rank 2 Drinfeld modules over $K$ with good reduction and $J_1 = j(\varphi^1)$, $J_2 = j(\varphi^2) \in W$; suppose
    $v(N) = 0$. For $k = 1, 2$ let
    \[
    n_k := \frac{1}{q-1}\#\Aut_{\bar\BF}(\bar \varphi^k) = \begin{cases}
        q + 1 & \text{if } v(J_k) > 0,\\
        1 & \text{if } v(J_k) = 0.
    \end{cases}
    \]
    Then the coefficients of $\Phi_N(X + J_1, Y+J_2) = \sum_{m,n}a_{m,n}X^mY^n\in W[X,Y]$ satisfy
    \[
    v(a_{m,n}) \geq n_1 C_{J_1,J_2}(N,\pi) - n_1 m - n_2 n = n_2 C_{J_2,J_1}(N,\pi) - n_1 m - n_2 n
    \]
    for all integers $m, n$.
\end{thm}

\begin{proof}
    Apply Proposition \ref{LSTG} to $\varphi^2$ with $M = 1$: this produces modules $\tilde{\varphi}^s$ with
    $v(j(\tilde{\varphi}^s) - J_2) = v(j(\tilde{\varphi}^s) - j(\tilde{\varphi}^l)) = \frac{1}{q-1}\#\Aut_{W/\pi}(\varphi^2) = n_2$, verifying hypothesis (2) of Proposition \ref{prop:general}. For hypothesis (3), let $\tilde\varphi$ be linked to $\tilde{\varphi}^s$ by a cyclic $N$-isogeny and suppose $v(j(\tilde\varphi) - J_1) > 0$. If $v(J_1) = 0$ then $n_1 = 1$ and there is nothing to prove. If $v(J_1) > 0$ then also $v(j(\tilde{\varphi})) > 0$; since $\tilde{\varphi}$ has good reduction (see Lemma \ref{lem:product}(ii)), Corollary \ref{cor:lattice} gives $v(j(\tilde{\varphi})) \geq q+1$, and likewise $v(J_1)\geq q+1$ since $J_1 = j(\varphi^1)$ with $\varphi^1$ of good reduction; hence $v(j(\tilde{\varphi}) - J_1) \geq q + 1 = n_1$. Proposition \ref{prop:general} now gives the first bound.

    For the second expression, it suffices to show $n_1C_{J_1,J_2}(N,\pi) = n_2C_{J_2,J_1}(N,\pi)$. Both sides equal $\frac{1}{q-1}$ times the number of cyclic $N$-isogenies $\bar \varphi^2 \to \bar \varphi^1$ over $\bar\BF$: indeed, the kernels counted by $C_{J_1,J_2}(N,\pi)$ are the cyclic $C\subset\bar \varphi^2[N]$ with $\bar \varphi^2/C\cong\bar \varphi^1$, and each such kernel is the kernel of exactly $\#\Aut_{\bar\BF}(\bar \varphi^1)$ isogenies $\bar \varphi^2\to\bar \varphi^1$; the analogous statement holds on the other side, and the two sets of isogenies are in bijection by Lemma \ref{lem:dual}.
\end{proof}

Theorem \ref{thm:main} follows by taking $J_1 = J_2 = 0$, $\varphi^1 = \varphi^2 = \phi^0$ and $W = \hat{A^{\mathrm{ur}}_P}$ for $\deg P$ odd: the coefficients of $\Phi_N(X,Y)$ lie in $A\subset W$, one has $n_1 = n_2 = q+1$, and $C_{0,0}(N,\pi) = C_0(N,P) = \psi(N)$ when $\deg P = 1$ by Proposition \ref{prop:fullpsi}. \qed

\begin{rem}\label{rem:vanishing}
    The argument does not use that $\deg P$ is odd: Theorem \ref{thm:main}(2) holds for every prime $P\nmid N$. For $\deg P$ even it carries no information, because $C_0(N,P) = C_0(N,0)$ by Proposition \ref{ordinary}(2), and all coefficients in its range vanish. More generally, let $J$ be algebraic over $k$ (for instance a singular modulus), write $\Phi_N(X+J,Y+J) = \sum_{m,n}a_{m,n}X^mY^n\in k(J)[X,Y]$ and $C_J(N,0) := \ord_X\Phi_N(X+J,J)$. Let $\mathfrak P$ be a prime of $k(J)$ above a monic prime of $A$, with $\mathfrak P\nmid N$ and $v_{\mathfrak P}(J) = 0$ (or $J = 0$); this excludes only finitely many primes. Let $W$ be the completion of the maximal unramified extension of the completion of $k(J)$ at $\mathfrak P$. By Corollary \ref{cor:lattice}, $J$ is the $j$-invariant of a Drinfeld module $\varphi^1 = \varphi^2$ over $W$ with good reduction, and Theorem \ref{thm:bigprimes} gives $v_{\mathfrak P}(a_{m,n})\geq n_1\big(C_{J,J}(N,\pi) - m - n\big)\geq C_J(N,0) - m - n$, since $n_1 = n_2\geq1$ and $C_{J,J}(N,\pi)\geq C_J(N,0)$. A non-zero element of $k(J)$ is divisible by only finitely many primes, so
    \[
    a_{m,n} = 0 \qquad\text{whenever } m + n < C_J(N,0),
    \]
    i.e.\ the curve $\Phi_N(X,Y) = 0$ has multiplicity exactly $C_J(N,0)$ at the point $(J,J)$. This explains the vanishing coefficients in Table \ref{tab:PhiT2T1} and in Section \ref{sec:computations}.
\end{rem}

\section{Ramified $j$-invariants and fractional exponents}\label{sec:ramified}

By Corollary \ref{cor:lattice}, a $j$-invariant $J$ with $0 < v_P(J) \not\equiv 0 \bmod (q+1)$ belongs to a Drinfeld module which acquires good reduction only over a ramified extension of $W$. In this section we prove a divisibility theorem for such $J$ at primes of degree one. Its constants are \emph{fractional} and exceed the value predicted by automorphisms alone; this is the Drinfeld analogue of the classical phenomenon that when good reduction requires a ramified base, the divisibility exponents acquire finer arithmetic structure (compare the exceptional rows of \cite[Table~1]{Breuer2025}, which are exactly the cases where good reduction requires a ramified base, with exponents such as $\tfrac{19}2, \tfrac92, \tfrac43, \tfrac12$, several of them depending on $N$ modulo $3$, $4$ or $6$). Remarkably, in the Drinfeld setting no complex multiplication hypothesis is needed.

By a change of variables we may restrict attention to the prime $P = T$: the $\BF_q$-algebra automorphism $\sigma : A\to A$, $T\mapsto T+\theta$, extends to an automorphism of $\mathbb C_\infty$ and transports rank $2$ Drinfeld modules, cyclic $N$-isogenies and $j$-invariants compatibly; consequently $\Phi_{\sigma(N)}$ is obtained from $\Phi_N$ by applying $\sigma$ to its coefficients, and the statement of Theorem \ref{thm:ramified} below for $(J, N)$ at the prime $T - \theta$ is equivalent to the statement for $(\sigma(J), \sigma(N))$ at the prime $T$.

\begin{thm}\label{thm:ramified}
    Let $J \in A$ with
    \[
    0 < v := v_T(J) < q+1, \qquad e := \frac{q+1}{\gcd(q+1,\, v)} \;(\geq 2),
    \]
    and let $N$ be monic with $T\nmid N$. Write $\Phi_N(X + J, Y + J) = \sum_{m,n}a_{m,n}X^mY^n \in A[X,Y]$. Then for all $m + n < \psi(N)$,
    \[
    v_T(a_{m,n}) \;\geq\; \Big(v + \frac{1}{e}\Big)\big(\psi(N) - m - n\big).
    \]
\end{thm}

Note two features of this statement. First, $C_J(N,T) = \psi(N)$ automatically, since $J\equiv0\bmod T$ and $j = 0$ is the only supersingular invariant in characteristic $T$; this is why $\psi(N)$ appears. Second, the constant exceeds by $\tfrac1e$ the value $v$ that the automorphism counts alone would give (via Proposition \ref{LSTG} over the ramified base); the excess is forced by the $\BF_{q^2}$-rationality of Lemma \ref{lem:Fq2}.

\begin{proof}
    Write $\pi'$ for a root of $x^e = \pi = T$ and $W' = W[\pi']$, a complete discrete valuation ring with the same (algebraically closed) residue field, totally tamely ramified of degree $e$ over $W$ (note $\gcd(e, q) = 1$ as $e \mid q+1$). Let $v'$ be the valuation of $W'$ with $v'(\pi') = 1$, so $v' = e\,v_T$ on $K$. By the Cohen structure theorem (see I.S. Cohen \cite{Coh46}), we may identify $W' = \bar\BF[[\pi']]$ and speak of the coefficients $x_i\in\bar\BF$ of $x = \sum_i x_i\pi'^i \in W'$.

    Put $z := ev/(q+1) \in \BZ$, so that $v'(J) = ev = z(q+1)$ and $1 \leq z = e\cdot\frac{v}{q+1} < e$, i.e.
    \begin{equation}\label{eq:se}
        z + 1 \;\leq\; e.
    \end{equation}
    By Corollary \ref{cor:lattice} over $W'$ there is a Drinfeld module $\rho$ over $W'$ with good reduction and $j(\rho) = J$; normalize $\rho_T = \gamma(T) + g_\rho\tau + \tau^2$, so $J = g_\rho^{q+1}$ and $v'(g_\rho) = z$. Let $\beta := (g_\rho)_z \neq 0$; then the lowest-order coefficient of $J$ is $(J)_{z(q+1)} = \beta^{q+1}$.

    \emph{Deformations.} Consider the modules
    \[
    \tilde{\varphi}^s : T\mapsto \gamma(T) + g_s\tau + \tau^2, \qquad g_s := g_\rho + \eps_s\pi'^{\,z+1},
    \]
    with $\eps_s\in \bar\BF^*$ pairwise distinct (via $W' = \bar\BF[[\pi']]$, $\bar\BF\subset W'$). Each $\tilde{\varphi}^s$ has good reduction and $\tilde{\varphi}^s\equiv\rho\bmod\pi'^{\,z+1}$. Since all these modules have $\Delta = 1$, we have $j(\tilde{\varphi}^s) = g_s^{q+1}$. As $q+1$ is prime to $p$, the group $\mu_{q+1}$ of $(q+1)$st roots of unity lies in $\bar\BF^*\subset W'^*$, and the factorization $x^{q+1} - y^{q+1} = \prod_{\zeta\in\mu_{q+1}}(x - \zeta y)$ gives
    \[
    j(\tilde{\varphi}^s) - j(\tilde{\varphi}^l) = \prod_{\zeta\in\mu_{q+1}}\big(g_s - \zeta g_l\big) = \prod_{\zeta\in\mu_{q+1}}\big((1-\zeta)g_\rho + (\eps_s - \zeta\eps_l)\pi'^{\,z+1}\big).
    \]
    For $\zeta = 1$ the factor is $(\eps_s - \eps_l)\pi'^{\,z+1}$, of valuation exactly $z+1$, since $\eps_s - \eps_l\in\bar\BF^*$. For $\zeta\neq1$ the factor is congruent to $(1-\zeta)g_\rho$ modulo $\pi'^{\,z+1}$. It therefore has valuation exactly $z$, whatever the $\eps$'s: its coefficient of $\pi'^{\,z}$ is $(1-\zeta)\beta\neq0$, because $1-\zeta\in\bar\BF^*$. Hence
    \[
    v'\big(j(\tilde{\varphi}^s) - j(\tilde{\varphi}^l)\big) = (z+1) + q\cdot z = z(q+1) + 1 \qquad (s\neq l).
    \]
    The same computation with $g_l$ replaced by $g_\rho$ (i.e.\ $\eps_l$ replaced by $0$) gives $v'(j(\tilde{\varphi}^s) - J) = z(q+1)+1$. So hypothesis (2) of Proposition \ref{prop:general} holds with $n_2 = z(q+1)+1$.

    \emph{Hypothesis (3) with $n_1 = z(q+1)+1$.} Let $C\subset\tilde{\varphi}^s[N]$ be a cyclic kernel with monic kernel polynomial $u = \sum u_i\tau^i$ and let $\tilde{\varphi} = \tilde{\varphi}^s/C$, with coefficient $\tilde g$ and $\tilde\Delta = 1$. All the modules $\tilde{\varphi}^s$ reduce to $\tau^2$ modulo $\pi'$ (as $v'(\gamma(T)) = e > 0$ and $v'(g_\rho) = z > 0$), so Lemma \ref{lem:Fq2} applies to $C$: $\bar u_0 \in\BF_{q^2}^*$, which gives $\bar u_0^{\,q^2-1} = 1$. By \eqref{eq:gprime-bottom},
    \[
    \tilde g = u_0^{1-q}\big(g_\rho + \eps_s\pi'^{\,z+1}\big) + u_0^{-q}u_1\big(T^q - T\big),
    \]
    in which $v'(T^q - T) = v'(T) + v'(T^{q-1}-1) = e \geq z+1$ by \eqref{eq:se}. Hence the coefficients of $\tilde g$ in degrees $\leq z$ come from the term $u_0^{1-q}g_\rho$ alone:
    \[
    v'(\tilde g) = z, \qquad (\tilde g)_z = \bar u_0^{1-q}\beta.
    \]
    Consequently $v'(j(\tilde{\varphi})) = v'(\tilde g^{\,q+1}) = z(q+1)$ with lowest-order coefficient
    \[
    \big(j(\tilde{\varphi})\big)_{z(q+1)} = \big((\tilde g)_z\big)^{q+1} = \bar u_0^{(1-q)(q+1)}\beta^{q+1} = \beta^{q+1} = (J)_{z(q+1)}.
    \]
    Therefore $v'\big(j(\tilde{\varphi}) - J\big) \geq z(q+1) + 1$ for \emph{every} cyclic kernel $C$ and every $s$ --- the hypothesis of the implication in Proposition \ref{prop:general}(3) is even superfluous. Since moreover every module cyclically $N$-isogenous to some $\tilde{\varphi}^s$ has supersingular reduction and hence $j \equiv 0 \equiv J\bmod\pi'$ (by the argument in the proof of Proposition \ref{prop:fullpsi}: $j = 0$ is the unique supersingular invariant in $A$-characteristic $T$), we have $C_{J,J}(N,\pi') = \psi(N)$, and Proposition \ref{prop:general} yields
    \[
    v'(a_{m,n})\;\geq\;\big(z(q+1)+1\big)\big(\psi(N)-m-n\big).
    \]
    As $v' = e\,v_T$ on $A$ and $\frac{z(q+1)+1}{e} = v + \frac1e$, the theorem follows.
\end{proof}

\begin{rem}
    The hypothesis $v_T(J) < q+1$ is needed: it guarantees \eqref{eq:se}, which keeps the term $u_0^{-q}u_1(T^q-T)$ away from the leading coefficient of $\tilde g$. For $v_T(J) \equiv 0\bmod(q+1)$ one is in the unramified situation of Theorem \ref{thm:bigprimes}, where the analogous interference does occur and the (optimal) constant is $q+1$. For comparison, the weaker bound $v_T(a_{m,n})\geq v\,(\psi(N)-m-n)$ follows from Theorem \ref{thm:main}(1) alone, by expanding $\Phi_N(X+J,Y+J) = \sum_{m',n'}a'_{m',n'}(X+J)^{m'}(Y+J)^{n'}$: the coefficient of $X^mY^n$ is a sum of terms of valuation $\geq(q+1)\max\big(0,\psi(N)-m-n-k\big) + vk$ with $k = m'+n'-m-n\geq0$, and since $v < q+1$ the minimum over $k$ is $v(\psi(N)-m-n)$. The content of Theorem \ref{thm:ramified} is thus the gain of $\frac1e$.
\end{rem}

\begin{corollary}\label{cor:isogenyinvariance}
    Let $P$ be a prime of degree one, $P\nmid N$, and let $J_1, J_2$ be algebraic over $k$ with $\Phi_N(J_1, J_2) = 0$; extend $v_P$ to $k(J_1,J_2)$. If $0\leq v_P(J_1) < q+1$, then $v_P(J_2) = v_P(J_1)$. If moreover $J_1\in A$ and $v := v_P(J_1) > 0$, then
    \[
    v_P(J_2 - J_1)\;\geq\; v + \frac1e, \qquad e = \frac{q+1}{\gcd(q+1, v)};
    \]
    in particular $J_1$ and $J_2$ have the same lowest-order term in their $P$-adic expansions.
\end{corollary}

\begin{proof}
    Write $\Phi_N = \sum a_{m,n}X^mY^n$ and $v_i = v_P(J_i)$. As $\Phi_N(J_1,Y)$ is monic in $Y$ with $P$-integral coefficients, $v_2\geq0$. Suppose $v_1\neq v_2$ and let $x\in\{J_1,J_2\}$ be the one with the smaller valuation $v_x$, $y$ the other one; then $v_x < q+1$, and $\Phi_N(x,y) = 0$ by symmetry of $\Phi_N$. The term $x^{\psi(N)}$ ($a_{\psi(N),0} = 1$) has valuation $\psi(N)v_x$. Every other term $a_{m,n}x^my^n$ has larger valuation: if $m+n < \psi(N)$, by Theorem \ref{thm:main}(1) its valuation is at least $(q+1)(\psi(N)-m-n) + mv_x + nv_y > \psi(N)v_x$, and if $m+n\geq\psi(N)$ then $n\geq1$ and its valuation is at least $mv_x+nv_y > \psi(N)v_x$. This contradicts $\Phi_N(x,y)=0$; hence $v_1 = v_2$.
    For the second claim we may take $P = T$ (by the change of variables $\sigma$ above). Write $\Phi_N(X + J_1, J_1) = \sum_m b_mX^m$, so that $b_m$ is the coefficient of $X^mY^0$ in $\Phi_N(X+J_1,Y+J_1)$. The roots of this polynomial are the differences $J_2 - J_1$; it is monic of degree $\psi(N)$, and $v_T(b_m)\geq(v+\frac1e)(\psi(N)-m)$ for $m<\psi(N)$ by Theorem \ref{thm:ramified}, so its Newton polygon shows that every root has valuation $\geq v+\frac1e$.
\end{proof}

\begin{rem}\label{rem:wang}
    Corollary \ref{cor:isogenyinvariance} is the Drinfeld analogue of Wang's application of the divisibility of classical modular polynomials to reduction types. For elliptic curves over a $2$-adic field linked by an isogeny of odd degree, with $0 < v_2(j(E_1)) < 15$, the argument above with $15$ in place of $q+1$ (the input being \cite[Thm.~1.1(1)]{Breuer2025}) gives $v_2(j(E_1)) = v_2(j(E_2))$ \cite[Prop.~2.18]{Wang2026}, and hence, in the case of potentially good supersingular reduction, $E_1$ and $E_2$ have the same Kodaira type \cite[\S3]{Wang2025}. In the Drinfeld setting the analogue of the Kodaira type of a module with potentially good reduction at $P$ is the order $e = (q+1)/\gcd(q+1,v_P(J))$ of the tame inertia character with values in $\Aut(\bar\varphi)=\BF_{q^2}^*$; by Corollary \ref{cor:isogenyinvariance}, as long as $v_P(J) < q+1$, both $e$ and the leading $P$-adic coefficient of $J$ are preserved by cyclic isogenies of degree prime to $P$, including those defined only over $\bar k$. The threshold $q+1$ is sharp: by \eqref{eq:gprime-bottom} with $u_1 = 1$, a module with $j = T^{q+2}$ (so $e = q+1$) is $(T-1)$-isogenous to modules whose $j$-invariants have $v_T = q+1$ (so $e = 1$).
\end{rem}

We now describe two families of \emph{singular moduli} to illustrate Theorem \ref{thm:ramified}; for these, the constant $v + \frac1e$ is exactly attained in all our computations (Section \ref{sec:comp-ramified}).

\subsection{CM by $\BF_q[\sqrt T\,]$}\label{sec:sqrtT}
Assume $q$ is odd and consider the imaginary quadratic extension $F = k(s)$, $s^2 = T$, with ring of integers $\CO_F = \BF_q[s]$, ramified precisely at $P = T$ (and at $\infty$). Since $\BF_q[s]$ is a principal ideal domain, there is up to isomorphism a unique Drinfeld module with CM by $\CO_F$, namely the rank 2 restriction of the rank 1 module $\rho_s = s + \tau$:
\[
\rho_T = (s + \tau)^2 = T + (s + s^q)\tau + \tau^2, \qquad
J_s := j(\rho) = (s + s^q)^{q+1} = T^{\frac{q+1}{2}}\big(1 + T^{\frac{q-1}{2}}\big)^{q+1} \in A.
\]
Here $v_T(J_s) = \frac{q+1}{2}$, so Theorem \ref{thm:ramified} applies with $v = \frac{q+1}2$, $e = 2$:
\[
v_T(a_{m,n}) \;\geq\; \frac{q+2}{2}\big(\psi(N) - m - n\big), \qquad T\nmid N.
\]
The supersingular primes of degree one for $J_s$ are $P = T$ (ramified in $F$) and the $P = T - \theta$ with $\theta$ a non-square in $\BF_q^*$ (inert in $F$). At the latter, $v_P(J_s) = q+1$ (each root of $1 + T^{(q-1)/2}$ is simple, raised to the power $q+1$), so Theorem \ref{thm:bigprimes} applies with $n_1 = n_2 = q+1$ and $C_J(N,P) = \psi(N)$. At split primes ($\theta$ a nonzero square) the reduction is ordinary and no divisibility occurs (Proposition \ref{ordinary}).

\subsection{CM by the conductor-$T$ order in $\BF_{q^2}[T]$}\label{sec:conductor}
The second family exists for \emph{all} $q$ and is the analogue of the classical singular moduli with non-maximal CM orders ($D = -12, -16, -27$ in \cite[Table~1]{Breuer2025}).

\begin{prop}\label{prop:conductor}
    Let $J_c := (-1)^q\,T\,(T^{q-1}-1)^{q+1} \in A$. Then:
    \begin{enumerate}
        \item $\Phi_T(X, 0) = (X - J_c)^{q+1}$;
        \item $J_c$ is the $j$-invariant of the Drinfeld modules with complex multiplication by the order $\CO' = A + T\,\BF_{q^2}[T]$ of conductor $T$ in $H = \BF_{q^2}(T)$; in particular $\CO'$ has class number one;
        \item $v_T(J_c) = 1$ and $v_{T-\theta}(J_c) = q+1$ for every $\theta\in\BF_q^*$.
    \end{enumerate}
\end{prop}

\begin{proof}
    Work over an algebraic closure of $k$. The $T$-torsion of $\phi^0$ is $$\phi^0[T] = \{x : Tx + x^{q^2} = 0\} = \BF_{q^2}\cdot\lambda,$$ where $\lambda^{q^2-1} = -T$. Its cyclic submodules (i.e. one-dimensional $\BF_q$-submodules) are the $q+1$ lines $C_\mu = \BF_q\,\mu\lambda$, $\mu\in\BF_{q^2}^*/\BF_q^*$, and $\Aut(\phi^0) = \BF_{q^2}^*$ permutes them transitively; hence the quotients $\phi^0/C_\mu$ are pairwise isomorphic and, by Lemma \ref{lem:product}, $\Phi_T(X,0) = (X - j(\phi^0/C))^{q+1}$ for any one of them. We compute $j(\phi^0/C)$ for $C = \BF_q\lambda$: the kernel polynomial is $u = \prod_{a\in\BF_q}(x - a\lambda) = x^q - \lambda^{q-1}x$, so $u_0 = -\lambda^{q-1}$, $u_1 = 1$, and by \eqref{eq:gprime-bottom} (with $g = 0$, $\Delta = 1$, whence also $\Delta' = 1$):
    \[
    j(\phi^0/C) = \big(u_0^{-q}(T^q - T)\big)^{q+1} = \frac{(T^q-T)^{q+1}}{(-1)^{q(q+1)}\lambda^{q(q^2-1)}},
    \]
    and using $\lambda^{q(q^2-1)} = (-T)^q$, we obtain
    \[
    j(\phi^0/C)
    = \frac{T^{q+1}(T^{q-1}-1)^{q+1}}{(-T)^{q}} = (-1)^q\,T\,(T^{q-1}-1)^{q+1}.
    \]
    This proves (1), and (3) follows by inspection.

    For (2): $\End(\phi^0/C)$ is an $A$-order in $\End^\circ(\phi^0) = H$ containing the image of $A + T\,\BF_{q^2}[T]$ (elements of $T\BF_{q^2}[T]$ kill $\phi^0[T]\supseteq C$, hence descend to the quotient), so $\End(\phi^0/C)\in\{\CO', \BF_{q^2}[T]\}$. If it were the maximal order then, $\BF_{q^2}[T]$ having class number one, $\phi^0/C\cong\phi^0$ and $J_c = 0$, which is false. So $\End(\phi^0/C) = \CO'$. Finally, every Drinfeld module $\rho$ with CM by $\CO'$ is isomorphic to some $\phi^0/C_\mu$. Indeed, over $\mathbb C_\infty$ we have $\rho\cong\phi^{\mathfrak a}$ for a fractional $\CO'$-ideal $\mathfrak a$ whose multiplicator ring is $\CO'$ (as in the proof of Lemma \ref{lem:CMtate}); put $\mathfrak a' := \BF_{q^2}[T]\,\mathfrak a \supseteq \mathfrak a$. The inclusion induces an isogeny $u : \phi^{\mathfrak a} \to \phi^{\mathfrak a'}\cong\phi^0$ (the multiplicator ring of $\mathfrak a'$ contains the maximal order $\BF_{q^2}[T]$, which has class number one) with kernel $e_{\mathfrak a}(\mathfrak a')\cong\mathfrak a'/\mathfrak a$, killed by $T$ because $T\,\BF_{q^2}[T]\subseteq\CO'$. So we have three cases: $\ker u$ is either trivial, or is all of $\phi^{\mathfrak a}[T]$, or is cyclic of degree $T$. In the first case $\BF_{q^2}[T]\,\mathfrak a = \mathfrak a$, so $\BF_{q^2}[T]$ would lie in the multiplicator ring $\CO'$ of $\mathfrak a$ --- this is not true. In the second, $\mathfrak a' = T^{-1}\mathfrak a$ would have the same multiplicator ring $\CO'$ as $\mathfrak a$, yet it contains $\BF_{q^2}[T]$ in its multiplicator ring --- this is not true. Hence $u$ must be in the third case, so $u$ is a cyclic $T$-isogeny, and its dual $\hat u : \phi^0\to\rho$ (Lemma \ref{lem:dual}) has cyclic kernel $\cong\rho[T]/\ker u$ of degree $T$; thus $\rho\cong\phi^0/C$ for a cyclic $T$-kernel $C\subset\phi^0[T]$, i.e.\ $\rho\cong\phi^0/C_\mu$ for some $\mu$. These quotients form a single isomorphism class; hence $h(\CO') = 1$.
\end{proof}

Since $v_T(J_c) = 1$, Theorem \ref{thm:ramified} applies with $v = 1$ and $e = q + 1$:

\begin{corollary}\label{cor:conductor}
    For monic $N$ with $T\nmid N$, the coefficients of $\Phi_N(X + J_c, Y+J_c)$ satisfy
    \[
    v_T(a_{m,n}) \;\geq\; \frac{q+2}{q+1}\,\big(\psi(N) - m - n\big), \qquad m + n < \psi(N),
    \]
    while at the primes $P = T - \theta$, $\theta\in\BF_q^*$, with $P\nmid N$, we have
    $v_P(a_{m,n}) \geq (q+1)(\psi(N) - m - n)$.
\end{corollary}

\begin{proof}
    The first bound is Theorem \ref{thm:ramified} with $v = 1$, $e = q+1$. For the second, $v_P(J_c) = q+1$ implies good reduction over the unramified $W = \hat{A_P^{\mathrm{ur}}}$ with $\bar J_c = 0$ (Corollary \ref{cor:lattice}), so Theorem \ref{thm:bigprimes} applies with $n_1 = n_2 = q+1$, and $C_{J_c,J_c}(N,\pi) = \psi(N)$ exactly as in Proposition \ref{prop:fullpsi}.
\end{proof}

For $q = 2$ the constant $\frac{q+2}{q+1} = \frac43$ coincides, amusingly, with the classical exponent $n_3 = \frac43$ proved for the discriminant $D=-27$ at $p = 3$ when $N\equiv\pm1\bmod 6$ \cite[Table~1]{Breuer2025}, which likewise belongs to a conductor-$p$ CM order whose singular modulus requires a ramified base for good reduction. The analogy is imperfect: the classical exponent depends on $N\bmod 6$ (it is $\frac12$ for $N\equiv\pm2\bmod 6$), and computations suggest that it is not optimal there (the expected optimal values being $\frac32$ and $1$, see \cite[\S1]{Breuer2025}), whereas $\frac{q+2}{q+1}$ is attained for every $N$ in our computations (Section \ref{sec:comp-ramified}).

\section{Computations}\label{sec:computations}

We have computed the Drinfeld modular polynomials $\Phi_N(X,Y)$ exactly, as elements of $A[X,Y]$, for
\[
\begin{aligned}
q = 2:&\quad N\in\big\{T,\ T+1,\ T^2,\ T(T+1),\ T^2+T+1,\ T^3,\ T^3+T+1\big\},\\
q = 3:&\quad N\in\big\{T,\ T+1,\ T+2,\ T(T+1),\ (T+1)^2,\ (T+2)^2,\ T^2+1,\\
&\phantom{\quad N\in\big\{}T^2+2,\ T^2+T+2,\ T^2+2T+2\big\},\\
q = 4,\,5,\,9:&\quad N = T,
\end{aligned}
\]
together with the univariate specializations $\Phi_N(X, 0)$ for all monic $N$ of degree $3$ and $4$ over $\BF_2$. All these computations use the proven height bound of Section \ref{sec:algorithm}, so every polynomial in our data set is certified to be exact; all numerical statements in this section refer to this data set. The implementation --- a Python library with verification scripts, an independent reimplementation in SageMath, and a companion Jupyter notebook that recomputes the polynomials from scratch and displays them in factored form --- is available at \cite{GitHubRepo}; together these verify every numerical claim of this section. The interested reader who would like to go further in the computational direction is referred to \cite{ABCLNP26}.

\subsection{The algorithm}\label{sec:algorithm}

Explicit results on Drinfeld modular polynomials in the literature include the closed form at level $T$ conjectured by Schweizer \cite{Schweizer1995} and proved by Bassa--Beelen \cite{BB2012} (see also \cite{BRaz2024}), and the algorithms of Caranay--Greenberg--Scheidler \cite{CGS2020}, who compute $\Phi_\ell$ modulo a prime, for $\ell$ irreducible. We compute $\Phi_N$ instead by \emph{deformation at a place of $A$}, running the proof of Theorem \ref{thm:main} at finite precision: the interpolation nodes are the congruence family of Proposition \ref{LSTG}, so the numerical verification below is performed by the very family of modules on which the proofs rest. The method is a function-field analogue of the deformation algorithm of Kunzweiler--Robert \cite{KR2024} for classical modular polynomials, with one simplification: since $A = \BF_q[T]$, a coefficient $a\in A$ is recovered \emph{exactly} from finitely many digits of its $P$-adic expansion $a = \sum_k r_kP^k$ (with $r_k\in A$, $\deg r_k < \deg P$) at a single place, so no Chinese remainder reconstruction over several places is required.

Fix a monic prime $P\nmid N$ of degree $d_P$ --- arbitrary: levels $N$ divisible by every prime of degree one, such as $N = T(T+1)$ for $q = 2$, admit no base place of degree one --- and let $W = \hat{A^{\mathrm{ur}}_P}$, with uniformizer $\pi = P$, valuation $v = v_P$ and residue field $\bar\BF$; write $\theta_P\in\BF_P = A/P\subseteq\bar\BF$ for the residue of $T$. The reduction of $\phi^0$ is $\bar\phi^0 : T\mapsto\theta_P+\tau^2$ over $\bar\BF$, and its $N$-torsion $\bar\phi^0[N]$, an $\BF_q$-space of dimension $2\deg N$, lies in a finite subfield of $\bar\BF$. Multiplication by $\theta_P$ commutes with $\tau^2$ only when $d_P\mid2$ (as $\tau^2\theta_P = \theta_P^{q^2}\tau^2$), so for $d_P\geq3$ the additive polynomial $\bar\phi^0_N$ must be evaluated as $N(\theta_P+\tau^2)$, a polynomial in the single $\BF_q$-linear operator $\theta_P+\tau^2$ (say by a Horner scheme), and not expanded coefficientwise as if $\theta_P$ were a scalar.

All computations take place in the finite ring $R = E[\pi]/(\pi^K)$, where $E = \BF_{q^m}$ is a large enough finite field. The ring $R$ receives its $A$-algebra structure from $\gamma(T) := t$, where $t$ is the unique root of the equation $P(t) = \pi$ with $t\equiv\theta_P\bmod\pi$ --- it exists by Hensel's lemma, $\theta_P$ being a simple root of the reduced equation, and Proposition \ref{prop:algcorrect} identifies $(R, \gamma)$ with a truncation of a subring of $W$. In particular $\gamma(P) = \pi$; and for $d_P = 1$, $P = T-\theta$, simply $t = \theta+\pi$. Kernel polynomials are assembled by the subspace recursion: writing $e_V(x) := \prod_{c\in V}(x-c)$ for a finite $\BF_q$-subspace $V$ of a domain, one has, for $v\notin V$,
\begin{equation}\label{eq:subspace}
    e_{V + \BF_q v} \;=\; e_V^q - e_V(v)^{q-1}\,e_V,
\end{equation}
since both sides are monic $\BF_q$-linear polynomials of degree $q^{\dim V + 1}$ vanishing on the $q^{\dim V+1}$ distinct points of $V + \BF_q v$. The full procedure is Algorithm \ref{alg:phiN}; note that its only divisions are by the units $P'(x)$ and $N(t)$ and by the Lagrange denominators, whose valuations are controlled exactly.

\begin{algorithm}[t]
\caption{$\Phi_N(X,Y)$ by deformation at a place of $A$}\label{alg:phiN}
\begin{algorithmic}[1]
\Require a prime power $q$; monic $N\in A$ of degree $d\geq1$; a monic prime $P\nmid N$ of degree $d_P$; an integer $h\geq h(\Phi_N)$
\Ensure the exact $\Phi_N(X,Y) = \sum_{m,n}a_{m,n}X^mY^n\in A[X,Y]$
\State $\psi\gets\psi(N)$; \quad $K\gets\lfloor h/d_P\rfloor + (q+1)\psi + 1$
\State choose $m\geq1$ with $d_P\mid m$, $\bar\phi^0[N]\subseteq E:=\BF_{q^m}$ and $\frac{q^m-1}{q+1}>\psi$; \quad $R\gets E[\pi]/(\pi^K)$
\State pick a root $\theta_P\in E$ of $P$; \ $t\gets$ the root of $P(t) = \pi$ in $R$ with $t\equiv\theta_P\bmod\pi$, by Newton iteration $x\gets x - P'(x)^{-1}\big(P(x)-\pi\big)$ \Comment{$\gamma(T) := t$; for $d_P = 1$: $t = \theta_P+\pi$}
\State compute $\bar\phi^0[N] = \ker\big(N(\theta_P+\tau^2)\big|_E\big)$ by $\BF_q$-linear algebra; enumerate its $\psi$ cyclic kernels $\bar C = A\cdot\bar w\cong A/NA$
\State choose $\eps_0,\dots,\eps_\psi\in E^*$ with pairwise distinct $(q+1)$-st powers
\For{$s = 0,\dots,\psi$}
    \State $\tilde{\varphi}^s : T\mapsto t + g_s\tau + \tau^2$ over $R$, \ with $g_s := \eps_s\pi$; \quad $y_s\gets j(\tilde{\varphi}^s) = g_s^{q+1}$
    \State lift each $\bar w\in\bar\phi^0[N]$ to the root $w\in R$ of $\tilde{\varphi}^s_N$ by Newton iteration $x\gets x - N(t)^{-1}\,\tilde{\varphi}^s_N(x)$
    \For{each cyclic kernel $\bar C = A\cdot\bar w$}
        \State $C\gets\BF_q\text{-span of } w,\ \tilde{\varphi}^s_T(w),\ \dots,\ \tilde{\varphi}^s_{T^{d-1}}(w)$ \Comment{the lift of $\bar C$}
        \State $u_C = \sum_{i=0}^{d}u_i\tau^i\gets e_C$ \ via \eqref{eq:subspace} \Comment{kernel polynomial, monic}
        \State $g'_C\gets g_s^{q^d} + u_{d-1} - u_{d-1}^{q^2}$; \quad $j'_C\gets(g'_C)^{q+1}$ \Comment{by \eqref{eq:gprime-top}, $\Delta = \Delta' = 1$}
    \EndFor
    \State $B^{(s)}(X)\gets\prod_{C}\big(X - j'_C\big)\in R[X]$ \Comment{$= \Phi_N(X, y_s)$ by Lemma \ref{lem:product}}
\EndFor
\State $\sum_{m,n}a_{m,n}X^mY^n\gets\sum_{s=0}^{\psi}B^{(s)}(X)\prod_{l\neq s}\frac{Y-y_l}{y_s-y_l}$ \Comment{correct modulo $\pi^{K-(q+1)\psi}$}
\State for each $(m,n)$: compute the $P$-adic digits $r_0,\dots,r_{\lfloor h/d_P\rfloor}\in A$ of $a_{m,n}$ by successive division: $x_0 = a_{m,n}$, \ $r_k\gets$ the element of degree $<d_P$ with $r_k(\theta_P)\equiv x_k\bmod\pi$, \ $x_{k+1}\gets\big(x_k - r_k(t)\big)/\pi$
\State \Return $\Phi_N(X,Y)$ with $a_{m,n} = \sum_{k}r_kP^k\in A$
\end{algorithmic}
\end{algorithm}

\begin{prop}\label{prop:algcorrect}
    Let $h\geq h(\Phi_N)$ and $K > h/d_P + (q+1)\psi(N)$. Then Algorithm \ref{alg:phiN} outputs the exact Drinfeld modular polynomial $\Phi_N(X,Y)\in A[X,Y]$, for every admissible choice of $P$, $E$, $\theta_P$ and $\eps_0,\dots,\eps_{\psi(N)}$.
\end{prop}

\begin{proof}
    As in the proof of Theorem \ref{thm:ramified}, use the Cohen structure theorem \cite{Coh46} to identify $W = \bar\BF[[\pi]]$, and embed $E$ into $\bar\BF$ so that the chosen root $\theta_P$ of $P$ maps to the residue of $T$ (possible, as all roots of $P$ in $\bar\BF$ are conjugate over $\BF_q$). The image $\gamma(T)\in W$ of $T$ satisfies $P(\gamma(T)) = \gamma(P) = \pi$ and $\gamma(T)\equiv\theta_P\bmod\pi$; since $\theta_P$ is a simple root of the reduction of $P(x)-\pi$ (the prime $P$ is separable), Hensel's lemma shows that $\gamma(T)$ is the unique such element of $W$, and the Newton iteration of step 3, which stays in $E[[\pi]]$ and converges quadratically, computes it: after $\lceil\log_2K\rceil$ iterations, $t\equiv\gamma(T)\bmod\pi^K$. This identifies $R = E[\pi]/(\pi^K)$, with the $A$-algebra structure of step 3, with a truncation of the subring $E[[\pi]]$ of $W$.

    Each $\tilde{\varphi}^s$ is then a Drinfeld module over $E[[\pi]]$ with $\Delta = 1$, hence has good reduction $\bar\phi^0$; moreover $y_s - y_l = \pi^{q+1}\big(\eps_s^{q+1} - \eps_l^{q+1}\big)$ has valuation \emph{exactly} $q+1$ for $s\neq l$, by the choice of the $\eps$'s (this is the congruence family of Proposition \ref{LSTG} with $M = 1$).

    By Lemma \ref{lem:product}(ii), the $N$-torsion points of $\tilde{\varphi}^s$ lie in $W$ and reduce bijectively and $A$-equivariantly onto $\bar\phi^0[N]\subseteq E$. The $\BF_q$-linear polynomial $\tilde{\varphi}^s_N(x)$ has derivative equal to its lowest coefficient $\gamma(N) = N(t)\in W^*$, a unit because $P\nmid N$; hence the Newton iteration of step 8, started at $\bar w$, remains in $E[[\pi]]$ and converges quadratically to the unique torsion point lifting $\bar w$, so after $\lceil\log_2 K\rceil$ iterations it is correct modulo $\pi^K$. Steps 9--14 are division-free polynomial expressions in these points, so by \eqref{eq:subspace}, \eqref{eq:gprime-top} and Lemma \ref{lem:product},
    \[
    B^{(s)}(X) \equiv \prod_{C\subset\tilde{\varphi}^s[N]\ \mathrm{cyclic}}\big(X - j(\tilde{\varphi}^s/C)\big) = \Phi_N\big(X, y_s\big) \pmod{\pi^K}.
    \]
    Write $\Phi_N(X,Y) = \sum_m b_m(Y)X^m$ with $b_m\in A[Y]$ of degree $\leq\psi(N)$. Since the $y_s$ are distinct, $b_m(Y) = \sum_s b_m(y_s)\prod_{l\neq s}\frac{Y - y_l}{y_s - y_l}$; each denominator is a unit times $\pi^{(q+1)\psi(N)}$ while the numerators are integral, so evaluating this formula with the values $b_m(y_s)$ known modulo $\pi^K$ --- which is step 16 --- determines every coefficient $a_{m,n}$ modulo $\pi^J$, where $J := K - (q+1)\psi(N) > h/d_P$.

    Finally, every $a\in A$ has a unique $P$-adic expansion $a = \sum_{k\geq0}r_kP^k$ with $r_k\in A$ of degree $<d_P$, and if $\deg_Ta\leq h$ then $r_k = 0$ for $k > h/d_P$. Evaluation at $\theta_P$ is a bijection from $\{r\in A : \deg r<d_P\}$ onto $\BF_P$ (it is reduction modulo $P$ followed by $A/P\cong\BF_P$), and applying $\gamma$ gives $\gamma(a) = \sum_kr_k(t)\pi^k$ with $r_k(t)\equiv r_k(\theta_P)\bmod\pi$. The successive divisions of step 17, applied to $a_{m,n}$ known modulo $\pi^J$, therefore recover the digits $r_0,\dots,r_{\lfloor h/d_P\rfloor}$ of $a_{m,n}$ in turn: identifying the $k$-th digit uses the current remainder $x_k\equiv\sum_{i\geq k}r_i(t)\pi^{i-k}$ modulo $\pi$ only, and each division by $\pi$ costs one digit of precision, which $J > h/d_P$ accommodates. As the remaining digits of $a_{m,n}$ vanish, step 18 returns $a_{m,n}$ exactly.
\end{proof}

When $d_P = 1$, say $P = T-\theta$, everything simplifies: $t = \theta+\pi$ exactly, $\bar\phi^0_N$ may be expanded coefficientwise, and the digits $r_k = c_k\in\BF_q$ are constants, so that step 18 returns $a_{m,n} = \sum_{k\leq h}c_k(T-\theta)^k$.

In practice we take for $h$ the integer part of the bound \eqref{explicitbound} with $b_q(N)$ replaced by $B_q$,
\[
h := \Big\lfloor\frac{q^2-1}{2}\,\psi(N)\big(\deg N - 2\lambda_N + B_q\big)\Big\rfloor \;\geq\; h(\Phi_N)
\]
by \cite[Thm.~1.1]{BPR2024}, so that Proposition \ref{prop:algcorrect} applies and every output is proven to be exact. We run with $K = \lceil h/d_P\rceil + (q+1)\psi(N) + 12$, i.e.\ with at least $11$ digits of precision beyond the minimum required by Proposition \ref{prop:algcorrect}, and the program checks that these extra digits of every coefficient vanish, that every residue $x_k\bmod\pi$ lies in the subfield $\BF_P$ of $E$, and that $\deg_T a_{m,n}\leq h$ for all $(m,n)$; by Proposition \ref{prop:algcorrect}, these checks can only fail through an error in the implementation. The bound is generous: over our data set it exceeds the actual height by a factor between $2.3$ and $5.9$ (at level $T$ and $q = 2$, for instance, $h = 45$, while $h(\Phi_T) = q^2(q+1) = 12$). The cost is dominated by the $(\psi(N)+1)\cdot q^{2\deg N}$ Newton lifts and the kernel-polynomial recursions, all on truncated power series of length $K = O\big(h(\Phi_N)/d_P\big)$; in a straightforward Python implementation on a laptop, levels of degree one take between one second and five minutes, the levels of degree $2$ and $3$ for $q = 2$ between half a minute and thirteen minutes, those of degree $2$ for $q = 3$ up to about an hour, and the largest case, $\Phi_T$ for $q = 9$ (computed at the place $T+1$, with residue field $\BF_{3^8}$), a few hours.

Three remarks. First, nothing in the algorithm requires the reduction $\bar\phi^0$ to be supersingular: the equidistance of the nodes $y_s$ is pure algebra, and the algorithm works at ordinary base places just as well --- it is the theorems, not the computation, that single out the supersingular ones. Second, since $\Phi_N$ does not depend on the chosen place, recomputing it at a second place is a strong end-to-end consistency check in which every intermediate object changes; we did this for several $(q, N)$, at places of degree up to $3$, always with identical output ($\Phi_{T(T+1)}$ for $q = 2$, where both degree-one primes divide $N$, was computed at a place of degree $2$). Third, the shifted polynomials $\Phi_N(X + J, Y + J)$ studied in Section \ref{sec:comp-ramified} are obtained from the exact $\Phi_N$ by binomial expansion.

The computed polynomials pass the following independent checks: each $\Phi_N$ is monic and symmetric, and satisfies the Kronecker congruence $\Phi_N\equiv(X^{|N|} - Y)(X - Y^{|N|})\bmod N$ for prime $N$ \cite[Thm.~3.1]{Bae1992}; the observed heights are consistent with \eqref{explicitbound} and, at level $T$, equal $q^2(q+1)$ in agreement with \cite{BRaz2024}; the specializations of our $\Phi_{T+1}$ and $\Phi_{T+2}$ for $q = 3$ agree with the tables of Caranay--Greenberg--Scheidler \cite[Ex.~6.2]{CGS2020}; our $\Phi_T$ agrees verbatim with the closed form of Bassa--Beelen \cite{BB2012} for $q\in\{2,3,4,9\}$; $\Phi_{T(T+1)}$ agrees with the resultant $\mathrm{Res}_Z\big(\Phi_T(X,Z), \Phi_{T+1}(Z,Y)\big)$ (checked at $25$ random specializations over $\BF_{2^{12}}$); and the independent SageMath implementation reproduces every polynomial coefficient for coefficient. Moreover, $\Phi_T$ for $q = 4$ and $q = 9$ was computed at two different base places ($T+1$ and $T+y$, where $\BF_q = \BF_p(y)$), with identical results (for $q=9$ the residue fields involved are $\BF_{3^8}$ and $\BF_{3^{16}}$).

\subsection{An example: $\Phi_{T^2+T+1}$ for $q = 2$}

\begin{table}[t]
{\small
\begin{tabular}{l}
$a_{2, 0} = \mathbf{T^{9}(T + 1)^{9}} \cdot (T^{2} + T + 1)^{3}$  \\
$a_{1, 1} = \mathbf{T^{9}(T + 1)^{9}} \cdot T^{2}(T + 1)^{2}$  \\
$a_{3, 0} = \mathbf{T^{6}(T + 1)^{6}} \cdot (T^{2} + T + 1)^{2}$  \\
$a_{2, 1} = \mathbf{T^{6}(T + 1)^{6}} \cdot (T^{2} + T + 1)(T^{4} + T + 1)$  \\
$a_{4, 0} = \mathbf{T^{3}(T + 1)^{3}} \cdot (T^{2} + T + 1)$  \\
$a_{3, 1} = \mathbf{T^{3}(T + 1)^{3}} \cdot (T^{2} + T + 1)$  \\
$a_{2, 2} = \mathbf{T^{3}(T + 1)^{3}} \cdot T^{2}(T + 1)^{2}(T^{2} + T + 1)(T^{4} + T + 1)$  \\
$a_{5, 0} = 1$ \\
$a_{4, 1} = a_{4,2} = T^{2} + T + 1$ \\
$a_{3, 2} = T(T + 1)(T^{2} + T + 1)(T^{3} + T^{2} + 1)(T^{3} + T + 1)$ \\
$a_{3, 3} = T^{2}(T + 1)^{2}(T^{2} + T + 1)$ \\
$a_{4, 4} = 1$; \quad $a_{0,0} = a_{1,0} = a_{5,1} = a_{5,2} = a_{4,3} = a_{5,3} = a_{5,4} = a_{5,5} = 0$ \\
\end{tabular}
\caption{The coefficients of $\Phi_{T^2+T+1}(X,Y) = \sum_{m,n}a_{m,n}X^mY^n$ for $q = 2$ (listing $a_{m,n} = a_{n,m}$ once), with the factors predicted by Theorem \ref{thm:main}(1) in bold. Note $C_0(N,0) = 2$, accounting for the two vanishing coefficients with $m + n < 2$; the remaining vanishing coefficients lie in the range $m + n > \psi(N) = 5$.}
\label{tab:PhiT2T1}
}
\end{table}

Table \ref{tab:PhiT2T1} displays $\Phi_{T^2+T+1}(X,Y)$ for $q=2$ in fully factored form; it is the analogue of \cite[Table~2]{Breuer2025}. Both primes of degree one are supersingular for $j = 0$ here, so Theorem \ref{thm:main}(1) makes every coefficient with $m + n < \psi(N) = 5$ divisible by $\big(T(T+1)\big)^{3(5 - m - n)}$; the factorizations show that this exponent is attained exactly at both primes, for instance at the coefficient $a_{2,0} = a_{C_0(N,0),\,0}$. The two vanishing coefficients with $m+n<2$ reflect $C_0(N,0) = \ord_X\Phi_N(X,0) = 2$, the number of kernels stable under $\End(\phi^0) = \BF_4[T]$ (namely the kernels of the two primes of $\BF_4[T]$ above $N$, which is of even degree and hence split). It is noteworthy how smooth the remaining cofactors are, in contrast with the classical situation \cite[Table~2]{Breuer2025}, where the unpredicted factors are typically large primes.

\subsection{Optimality of the constant $q+1$}\label{sec:optimality}

For each computed $\Phi_N$ and each prime $P\nmid N$, define
\[
n_{\mathrm{opt}}(N, P) := \min\bigg\{\frac{v_P(a_{m,n})}{C_0(N,P) - m - n} \;:\; m+n < C_0(N,P),\ a_{m,n}\neq0\bigg\},
\]
the best constant that could replace $q+1$ in Theorem \ref{thm:main} for that pair $(N, P)$; recall that $C_0(N,P) = \psi(N)$ when $\deg P = 1$. The theorem asserts $n_{\mathrm{opt}}(N,P)\geq q+1$ for $\deg P$ odd. In our data this is an equality in every single case: for all $44$ pairs $(N, P)$ with $\deg P = 1$, $P\nmid N$, arising from the list above ($q\in\{2,3,4,5,9\}$), we found
\[
n_{\mathrm{opt}}(N, P) = q+1.
\]

The minimum is not always attained at the lowest coefficient. By Lemma \ref{lem:product}, the first nonvanishing coefficient of $\Phi_N(X, 0)$ is $a_{c,0} = \pm\prod_C j(\phi^0/C)$, with $c = C_0(N,0)$ and $C$ running over the cyclic kernels with $j(\phi^0/C)\neq0$; and for $\deg P = 1$, formulas \eqref{eq:gprime-bottom} and \eqref{eq:vTqT} (with $g = 0$, $\Delta = 1$) give
\[
v_P\big(j(\phi^0/C)\big) = (q+1)\big(1 + v_P(S_C)\big),
\]
with $S_C = u_1/u_0$ as in \eqref{SC}. Thus $n_{\mathrm{opt}}(N,P) = q+1$ is witnessed by $a_{c,0}$ exactly when $v_P(S_C) = 0$ for every such kernel. This is what happens in $42$ of our $44$ cases, but not in all: for $(q, N, P) = (2,\ T^3+T+1,\ T+1)$, three of the nine kernels satisfy $P\mid S_C$, whence $v_P(a_{0,0}) = 36 > 27 = (q+1)\psi(N)$, and for $(q,N,P) = (2,\ T^3,\ T+1)$ six of the twelve kernels do, whence $v_P(a_{0,0}) = 54 > 36$. In both cases the minimum $q+1$ is attained at other coefficients instead. The optimality of $q+1$ is thus a subtler phenomenon than the nonvanishing of the residues $\bar S_C$.




\subsection{Higher-degree primes}\label{sec:higherdeg}

For $\deg P = 2$ and $J = 0$ the reduction is ordinary, and indeed $C_0(N, P) = C_0(N,0)$ in all computed cases, with all coefficients $a_{m,n}$, $m + n < C_0(N,0)$, vanishing identically, in accordance with Proposition \ref{ordinary}(2) (the CM order $\BF_{q^2}[T]$ of $\phi^0$ is maximal) and Remark \ref{rem:vanishing}.

At primes of degree $3$ (supersingular for $j = 0$), Proposition \ref{supersingular}(2) constrains the excess $C_0(N,P) > C_0(N,0)$: the CM order $\BF_{q^2}[T]$ of $\phi^0$ has unit discriminant, so for odd $q$ the bound reads $\deg P\leq\deg N$. The exhaustive computation reported here has $q = 2$, which the proposition as proved does not cover (a characteristic-$2$ analogue of Lemma \ref{lem:kaneko} would require a separate treatment of the quadratic form); the $q=2$ data turn out to obey the same bound. To probe this we computed the univariate polynomials $\Phi_N(X, 0)\in A[X]$ for \emph{all} $24$ monic $N$ of degree $3$ and $4$ over $\BF_2$ (the six levels divisible by $T(T+1)$ at places of degree $2$ or $3$), and determined $C_0(N,P)$ at every prime $P\nmid N$ of degree $\leq 3$. At the primes of degree $3$ we also recomputed $C_0(N,P)$ independently, as the number of cyclic kernels $\bar C\subset\bar\phi^0[N]$ in characteristic $P$ with $j(\bar\phi^0/\bar C) = 0$, a computation which needs no height bound; the two counts agree. Genuine excess does occur, in perfect agreement with the degree bound: $20$ of the $24$ levels have excess at some prime of degree $3$, for instance
\[
\begin{aligned}
C_0(T^3,\, T^3+T^2+1) &= 3, &\qquad C_0(T^4,\, T^3+T+1) &= 6,\\
C_0(T^4+T^3+1,\, T^3+T^2+1) &= 5, &\qquad C_0(T^4+T^3,\, T^3+T^2+1) &= 9,
\end{aligned}
\]
while $C_0(N,0) = 0, 0, 2, 0$ respectively. The excess $C_0(N,P) - C_0(N,0)$ is always a multiple of $q+1$, as it must be: the kernels it counts are not stable under $\Aut(\phi^0) = \BF_{q^2}^*$, which permutes them in orbits of length $q+1$ and preserves $j(\phi^0/C)$. (Curiously, the four levels without excess at a prime of degree $3$ are the two \emph{prime} levels $N$ of degree $3$, the prime $T^4+T+1$, and $(T^2+T+1)^2$; for $N = T(T+1)$, too, $C_0(N,P) = 0$ at both primes of degree $3$.)
We verified Theorem \ref{thm:main}(2) directly on the full polynomial $\Phi_{T^3}(X,Y)$: at the excess prime $P = T^3+T^2+1$ its coefficients satisfy $v_P(a_{m,n})\geq(q+1)\big(C_0(N,P) - m - n\big) = 3(3 - m - n)$, and $n_{\mathrm{opt}}(T^3, P) = 3 = q+1$, so the constant is again optimal. This is the only setting in our data where part (2) of Theorem \ref{thm:main} is tested with $0 < C_0(N,P) < \psi(N)$, and it, too, is sharp.

\subsection{Ramified and singular moduli}\label{sec:comp-ramified}

For $q = 3$ and the singular modulus $J = J_s = T^2(1+T)^4$ of Section \ref{sec:sqrtT} we computed the shifted polynomials $\Phi_N(X + J, Y+J)$ for
\[
N\in\big\{T+1,\ T+2,\ T^2+1,\ T(T+1),\ (T+1)^2,\ (T+2)^2,\ T^2+2,\ T^2+T+2,\ T^2+2T+2\big\}
\]
and tabulated $n_{\mathrm{opt}}$ --- defined as in Section \ref{sec:optimality}, now for the coefficients of $\Phi_N(X+J, Y+J)$ and with $C_0(N,P)$ replaced by $C_J(N,P)$ --- at each prime of degree $\leq2$ not dividing $N$. The results are uniform across all nine values of $N$ (prime, composite, and non-squarefree alike; each prime $P$ is of course considered only when $P\nmid N$) and confirm the theoretical picture exactly --- in every case each bound is attained:

\vspace{0.5cm}

\begin{center}
{\small
\begin{tabular}{lllll}
\toprule
$P$ & splitting in $k(\sqrt T)$ & $v_P(J_s)$ & $C_J(N,P)$ & $n_{\mathrm{opt}}$ \\
\midrule
$T$ & ramified & $2 = \frac{q+1}{2}$ & $\psi(N)$ & $5/2 = \frac{q+2}{2}$ (Thm.\ \ref{thm:ramified})\\
$T+1$ & inert & $4 = q+1$ & $\psi(N)$ & $4 = q+1$ (Thm.\ \ref{thm:bigprimes})\\
$T+2$ & split & $0$ & $C_J(N,0)$ & --- (ordinary; Prop.\ \ref{ordinary})\\
$T^2+T+2$, $T^2+2T+2$ & (deg 2, $J\equiv j_P$) & $0$ & $\psi(N)$ & $1$ (Prop.\ \ref{prop:deg2})\\
$T^2 + 1$ & (deg 2, $J\not\equiv j_P$) & $0$ & $C_J(N,0)$ & --- (ordinary)\\
\bottomrule
\end{tabular}
}
\end{center}

\vspace{0.5cm}

In particular the constant $\frac{q+2}{2}$ of Theorem \ref{thm:ramified} is attained for every $N$ in the list: in contrast with the classical ramified cases (where, e.g., for $D = -4$ the exponent at $p = 2$ depends on $N\bmod4$ \cite[Table 1]{Breuer2025}), no dependence on $N$ is visible here.

We also tested Theorem \ref{thm:ramified} in its full generality, i.e.\ for $j$-invariants $J$ with $0 < v_T(J) < q+1$ which are \emph{not} singular moduli. For
\[
q = 2:\; J \in \{T,\ T^2,\ T(1+T),\ T + T^4\}, \qquad
q = 3:\; J \in \{T,\ T^2,\ T^3\},
\]
paired with every level $N$ above with $T\nmid N$ (three levels for $q=2$ and eight for $q=3$, so $36$ pairs $(J,N)$ in all), the bound $\big(v_T(J) + \tfrac1e\big)\big(\psi(N) - m - n\big)$ holds and the constant
\[
v_T(J) + \tfrac1e \;\in\; \big\{\tfrac43,\ \tfrac73\big\}\ (q=2), \qquad
\big\{\tfrac54,\ \tfrac52,\ \tfrac{13}4\big\}\ (q = 3)
\]
is \emph{exactly attained} in every case --- covering all nonzero residue classes of $v_T(J)$ modulo $q+1$.

For the conductor-$T$ singular modulus $J_c$ of Proposition \ref{prop:conductor} we tested all levels $N$ above (with $T\nmid N$) for $q\in\{2,3\}$ and found, again uniformly in $N$:
at $P = T$ the constant $\tfrac{q+2}{q+1}$ of Corollary \ref{cor:conductor} is exactly attained ($\tfrac43$ for $q=2$, $\tfrac54$ for $q=3$);
at $P = T-\theta$, $\theta\in\BF_q^*$, where $v_P(J_c) = q+1$, the constant $q+1$ is exactly attained;
and at primes of degree $2$ (which split in $H$, so that $J_c$ has ordinary reduction there) one finds $C_J(N,P) = C_J(N,0) \in\{0,2\}$ with all coefficients $a_{m,n}$, $m+n < C_J(N,0)$, vanishing identically, as required by Proposition \ref{ordinary}(2) (these primes do not divide the conductor $T$ of $\CO'$) and Remark \ref{rem:vanishing}.


\end{document}